\documentclass[draft,10pt]{article}
\usepackage{amssymb}
\usepackage{verbatim,cite}
\usepackage{latexsym}
\usepackage{amsmath}
\usepackage{eucal}

\newtheorem{sub}{\name}[section]
\newcommand{\bs}{
\begin{sub}}
\newcommand{\es}{
\end{sub}}
\newcommand{\bsl}[1]{
\begin{sub}\label{#1}}
\newcommand{\bth}[1]{\def\name{Theorem}
\begin{sub}\label{t:#1}}
\newcommand{\blemma}[1]{\def\name{Lemma}
\begin{sub}\label{l:#1}}
\newcommand{\bcor}[1]{\def\name{Corollary}
\begin{sub}\label{c:#1}}
\newcommand{\bdef}[1]{\def\name{Definition}
\begin{sub}\label{d:#1}}
\newcommand{\bprop}[1]{\def\name{Proposition}
\begin{sub}\label{p:#1}}
\newcommand{\brem}[1]{\def\name{Remark}
\begin{sub}\label{r:#1}}
\newcommand{\bex}[1]{\def\name{Example}
\begin{sub}\label{e:#1}}

\newcommand{\rth}[1]{Theorem~\ref{t:#1}}
\newcommand{\rlemma}[1]{Lemma~\ref{l:#1}}

\newcommand{\rprop}[1]{Proposition~\ref{p:#1}}


\newcommand{\BA}{
\begin{array}}
\newcommand{\EA}{
\end{array}}

\newcommand{\BAN}{\renewcommand{\arraystretch}{1.2}
\setlength{\arraycolsep}{2pt}
\begin{array}}
\newcommand{\BAV}[2]{\renewcommand{\arraystretch}{#1}
\setlength{\arraycolsep}{#2}
\begin{array}}
\newcommand{\BSA}{
\begin{subarray}}
\newcommand{\ESA}{
\end{subarray}}
\newcommand{\BAL}{
\begin{aligned}}
\newcommand{\EAL}{
\end{aligned}}
\newcommand{\BALG}{
\begin{alignat}}
\newcommand{\EALG}{
\end{alignat}}
\newcommand{\BALGN}{
\begin{alignat*}}
\newcommand{\EALGN}{
\end{alignat*}}
\newcommand{\note}[1]{\textit{#1.}\hspace{2mm}}
\newcommand{\Proof}{\note{Proof}}
\newcommand{\qeda}{\hspace{10mm}\hfill $\square$}

\newcommand{\Remark}{\note{Remark}}

\newcommand{\forevery}{\quad \forall}

\newcommand{\abs}[1]{\left |#1\right |}

\def\angb<#1>{\langle #1 \rangle}

\newcommand{\myfrac}[2]{{\displaystyle \frac{#1}{#2} }}
\newcommand{\myint}[2]{{\displaystyle \int_{#1}^{#2}}}
\newcommand {\dint}{{\displaystyle \int\!\!\int}}

\newcommand{\prt}{
\partial}

\newcommand{\ti}{\times}

\def\ga{\alpha}     \def\gb{\beta}       
       \def\gd{\delta}      \def\ge{\epsilon}
\def\gth{\theta}                         
\def\gf{\phi}           
      \def\gk{\kappa}      \def\gl{\lambda}
        \def\gn{\nu}         \def\gp{\pi}
    \def\gr{\rho}        
       \def\gt{\tau}
      \def\gw{\omega}
                \def\gz{\zeta}
\def\Gg{\Gamma}     \def\Gd{\Delta}      \def\Gf{\Phi}

\def\Gl{\Lambda}

      \def\CP{{\mathcal P}}
\def\CA{{\mathcal A}}      
\def\CD{{\mathcal D}}

      \def\CV{{\mathcal V}}

   \def\BBE {\mathbb E}

   \def\BBN {\mathbb N}    
   \def\BBR {\mathbb R}

\allowdisplaybreaks
\begin{document}
\title {\bf  Singular solutions of some nonlinear parabolic equations with spatially inhomogeneous absorption}
\author{{\bf\large Andrey Shishkov}
\hspace{2mm}\vspace{3mm}\\
{\it \normalsize  Institute of Applied Mathematics and Mechanics of NAS of Ukraine},\\
{\it\normalsize R. Luxemburg str. 74, 83114 Donetsk, Ukraine}\\
\vspace{3mm}\\
{\bf\large Laurent V\'eron}
\vspace{3mm}\\
{\it\normalsize Laboratoire de Math\'ematiques et Physique Th\'eorique},
\\
{\it\normalsize Universit\'e Fran\c{c}ois-Rabelais,  37200 Tours,
France}}

\date{}
\maketitle
\begin{center} {\bf\small Abstract}

\vspace{3mm}
\hspace{.05in}
\parbox{5in} {{\small  We study the limit behaviour of solutions of
$\prt_tu-\Gd u+h(\abs x)\abs u^{p-1}u=0\quad\text { in }\BBR^N\ti (0,T)
$ with initial data $k\gd _{0}$ when $k\to\infty$, where $h$ is a positive nondecreasing function and $p>1$.
If $h(r)=r^{\gb}$, ($\gb>-2$) we prove that the limit function $u_{\infty}$ is an explicit very singular solution.
If $\liminf_{r\to 0}r^2\ln (1/h(r))>0$, $u_{\infty}$ has a persistent singularity at $(0,t)$ ($t\geq 0$).
If $\int_{0}^{r_{0}}r\ln (1/h(r))\,dr<\infty$, $u_{\infty}$ has a pointwise singularity localized at $(0,0)$.}}
\end{center}

\section{Introduction}
Consider
\begin{equation}\label{E1}
\prt_tu-\Gd u+h( x)\abs u^{p-1}u=0\quad\text { in }Q_{T}:=\BBR^N\ti (0,T),
\end{equation}
with $p>1$ and $h$ is a nonnegative measurable function defined in $\BBR^N$. It is well known that if
\begin{equation}\label{int1}
\dint_{Q_{T}}h(x)E^p(x,t)dx\,dt<\infty,
\end{equation}
where $E(x,t)=(4\gp t)^{-N/2}e^{-|x|^2/4t}$ is the heat kernel,
then, for any $k>0$ there exists a unique solution $u=u_{k}$ to
(\ref{E1}) satisfying initial condition
\begin{equation}\label{init}
u(.,0)=k\gd_{0}
\end{equation}
in the sense of measures in $\BBR^N$. Furthermore the mapping
$k\mapsto u_{k}$ is increasing. If it assumed that $h$ is positive
essentially locally bounded from above and from below in
$\BBR^N\setminus\{0\}$, then the set  $\{u_{k}\}$ is also bounded
in the $C^1_{loc}(\overline
Q_{T}\setminus\{0\times(0,\infty)\})$-topology. Thus there exist
$u_{\infty}:=\lim_{k\to\infty}u_{k}$  and $u_{\infty}$ is a
solution of (\ref{E1}) in $Q_{T}\setminus\{0\ti (0,\infty)\}$.
 Furthermore $u_{\infty}$ is continuous in $\overline
Q_{T}\setminus\{0\ti [0,\infty)\}$ and vanishes on
$\BBR^N\setminus\{0\}\ti \{0\}$. We shall prove that only two
situations can occur: \smallskip

\noindent (i) Either $u_{\infty}(0,t)$ is finite for every $t>0$ and $u_{\infty}$ is a solution of (\ref{E1})
in $Q_{T}$. Such a solution which has a pointwise singularity at $(0,0)$ is called a {\it very singular solution}
(abr. V.S.S.)\smallskip

\noindent (ii) Or $u_{\infty}(0,t)=\infty$ for every $t>0$ and $u_{\infty}$ is a solution of (\ref{E1})
in $Q_{T}\setminus\{0\ti (0,\infty)\}$ only. Such a solution with a persistent singularity is called a
{\it razor blade} (abr. R. B.).\smallskip

\noindent In the well-known article \cite{BPT}, Brezis, Peletier
and Terman proved in 1985 that $u_{\infty}$ is a V.S.S., if
$h(x)\equiv 1$. Furthermore they showed that
$u_{\infty}(x,t)=t^{-1/(p-1)}f(x/\sqrt t)$ for $(x,t)\in Q_{T}$
where $f$ is the unique positive (and radial) solution of the
problem
\begin{equation}\label{int2}\left\{\BA {l}
-\Gd f -\myfrac{1}{2}\eta.\nabla f-\myfrac{1}{p-1}f+\abs f^{p-1}f=0\quad\text{in }\BBR^N
\\[2mm]
\phantom{-----.}
\lim_{|\eta|\to\infty}|\eta|^{2/(p-1)}f(\eta)=0.
\EA\right.\end{equation} Their proof of existence and uniqueness relied on
shooting method in ordinary differential equations (abr. O.D.E.). The already mentioned
self-similar very singular solutions of the problem
(\ref{int2}) was discovered independently in \cite{GKS} too.
 Later on, a new proof of existence, has been given by Escobedo and Kavian
\cite{EK} by a variational method in a weighted Sobolev space.
More precisely they proved that the following functional
\begin{equation}\label{int3}
v\mapsto J(v)=\myfrac{1}{2}\int_{\mathbb{R}^N}\left(|\nabla
v|^2-\myfrac{1}{p-1}v^2+\myfrac{2}{p+1}|v|^{p+1}\right)
K(\eta)d\eta
\end{equation}
achieves a nontrivial minimum in $H^1_{K}(\BBR^N)$, where $K(\eta)=e^{|\eta|^2/4}$.

In this article we first study equation (\ref{E1}) when
$h(x)=|x|^{\gb}$ ($\gb\in\BBR$). Looking for self-similar solutions
under the form $u(x,t)=t^{-(2+\gb)/2(p-1)}f(x/\sqrt t)$, we are
led to
\begin{equation}\label{int4}\left\{\BA {l}
-\Gd f-\myfrac{1}{2}\eta.\nabla f-\myfrac{2+\gb}{2(p-1)}f+|\eta|^\gb\abs f^{p-1}f=0\quad\text{in }\BBR^N\\
[4mm]
f\in H_{loc}^1(\BBR^N)\cap L_{loc}^{p+1}(\BBR^N;|\eta|^\gb d\eta)\cap C^2(\BBR^N\setminus\{0\})\\
[2mm]
\phantom{-----.}\lim_{|\eta|\to\infty}|\eta|^{(2+\beta)/(p-1)}f(\eta)=0,
\EA\right.\end{equation} and the associated functional
\begin{equation}\label{int5}
v\mapsto J(v)=\myfrac{1}{2}\int_{\mathbb{R}^N}\left(|\nabla
v|^2-\myfrac{2+\gb}{2(p-1)}v^2+\myfrac{2}{p+1}|\eta|^\gb|v|^{p+1}\right)
K(\eta)d\eta.
\end{equation}
We prove the following\medskip

\noindent {\bf Theorem A}{ \it I- Assume $\gb\leq N(p-1)-2$; then there exists no nonzero solution to (\ref{int4}).
\smallskip

\noindent II- Assume $\gb>N(p-1)-2$; then there exists a unique positive solution $f^*$ to (\ref{int4}).}\medskip

One of the key arguments in the study of isolated singularities of  (\ref{E1}) is the following a priori estimate
\begin{equation}
|u(x,t)|\leq \myfrac{\tilde c}{(t+|x|^2)^{(2+\gb)/2(p-1)}}\forevery (x,t)\in Q_{T}
\end{equation}
valid for any $p>1$ and $\gb>-2$. The remarkable aspect of this
proof is that it is based upon the auxiliary construction of the
maximal solution of (\ref{E1}) under a selfsimilar form. Next we
give two proofs of II, one based upon scaling transformations and
asymptotic analysis of O.D.E., combining ideas from \cite {BPT},
\cite{BO} and \cite{MV0}, and the second based on variational
methods, extending some ideas from \cite{EK} and valid in a more general context. As a consequence we
prove\medskip

\noindent {\bf Theorem B}{ \it Assume $\gb>N(p-1)-2$, then $u_{\infty}(x,t)=t^{-(2+\gb)/2(p-1)}f^*(x/\sqrt t)$.}\medskip

It must be noticed that, if $\gb\leq N(p-1)-2$, $u_{k}$ does not
exist, and more precisely, the isolated singularities of solutions
of (\ref{E1}) are removable.

Next we consider the case of more degenerate potentials $h(x)$:
\begin{equation}\label{Flat}
\frac{h(x)}{|x|^\alpha}\to 0 \ \ \ \text{as} \ \ \ |x|\to 0 \ \
\forall \alpha >0.
\end{equation}
In the set of such potentials we find the borderline which
separates the above mentioned two possibilities (i) --- (V.S.S.) and
(ii) --- (R.B). Remark that in the case of flat potentials like
\eqref{Flat}, the corresponding solution $u_\infty(x,t)$ does not have
self-similar structure and we haveto find some alternative
techniques for the study of the structure of $u_\infty$. Main
results of the paper are the following two statements.
\medskip

\noindent {\bf Theorem C (sufficient condition of V.S.S.
solution)}{ \it Assume that the function $h$ is continuous and
positive in $\mathbb{R}^N\setminus\{0\}$ and verifies the following
flatness condition
\begin{equation}\label{Fl1}
|x|^2 \ln \left(\frac{1}{h(x)}\right)\leq
\omega(|x|)\Leftrightarrow h(x) \geq e^{-\omega(|x|)/|x|^2} \ \
\forall x\in \mathbb{R}^N,
\end{equation}
where the function $\gw\geq 0$ is nondecreasing, satisfies the
following Dini-like condition
\begin{equation}\label{int8}
\int_{0}^1\myfrac{\gw(s) ds}{s}<\infty,
\end {equation}
and the additional technical condition
\begin{equation}\label{int9}
s\gw'(s)\leq (2-\ga_{0})\gw(s) \quad\text {near }0,
\end {equation}
for some $\ga_{0}\in (0,2)$. Then $u_{\infty}(x,t)<\infty$ for any
$(x,t)\in Q_{T}$. Furthermore there exists positive constants
$C_{i}$ ($i=1,2,3$), depending only on $N$, $\alpha_0$ and $p$,
such that
\begin{equation}\label{int10}
\int_{\BBR^N}u^2_{\infty}(x,t)\,dx\leq
C_{1}t\exp\left[C_{2}\left(\Gf^{-1}\left(C_{3}t\right)\right)^{-2}\right]\forevery
t>0,
\end {equation}
where $\Gf^{-1}$ is the inverse function of
$$\Phi(\gt) : = \int_{0}^\gt\myfrac{\gw(s)}{s}\,ds.
$$
} Notice that (\ref{int8})-(\ref{int9}) is satisfied if $h(x)\geq
Ce^{-\abs x^{\gth-2}}$ for some $\gth>0$.\medskip

\noindent {\bf Theorem D (sufficient condition of R.B. solution)}{
\it Assume $h$ is continuous and positive  in
$\BBR^N\setminus\{0\}$ and satisfies
\begin{equation}\label {int6}
 \liminf_{x\to 0}|x|^2\ln
\left(\frac{1}{h(x)}\right)>0\Leftrightarrow\exists \omega_0=\text{\rm const} >0 : h(x)\leq
\exp \left(-\frac{\omega_0}{|x|^2}\right).
\end{equation}
Then $u_{\infty} (0,t)=\infty$ for any $t>0$, and $t\mapsto u_{\infty} (x,t)$ is increasing.
If we denote $U(x)=\lim_{t\to\infty}u_{\infty} (x,t)$, then $U$ is the minimal large solution of
\begin{equation}\label {int7}
-\Gd u+h(x)u^p=0\quad\text{in }\BBR^N\setminus\{0\},
\end{equation}
i.e. the smallest solution of (\ref{int7}) which satisfies
\begin{equation}\label {int7*}
\myint{B_{\ge}}{}u(x)dx=\infty\quad\forall\ge>0.
\end{equation}
}

Theorem C is proved by some new version of local energy method. A
similar variant of this method was used in \cite{B-Sh} for the
study of extinction properties of solutions of nonstationary
diffusion-absorption equations.
\medskip

Theorem D is obtained by constructing local appropriate
sub-solutions. The monotonicity and the limit property of
$u_{\infty}$ are characteristic of razor blades solutions \cite
{VV}.
\medskip

A natural question which remains unsolved is to characterize $u_\infty$ if the potential
$h(x)$ satisfies
$$
h(x) \approx \exp \left(-\frac{\omega(|x|)}{|x|^2}\right),
$$
where $\omega(s)\to 0$ as $s\to 0$ and
$$\int_{0}^1\frac{\omega(s)ds}{s}=\infty.$$

This article is the natural continuation of \cite{MV2}, \cite{ShV} where (\ref{E1}) is replaced by
\begin{equation}\label{E1''}
\prt_tu-\Gd u+h( t)\abs u^{p-1}u=0\quad\text { in }Q_{T}.
\end{equation}
In equation (\ref{E1''}), the function $h\in C([0,T])$ is positive in
$(0,T]$ and vanishes only at $t=0$. In the particular case
$h(t)=t^\gb$ ($\gb>0$), $u_{k}$ exists if and only if
$1<p<1+2(1+\gb)/N$,  and $u_{\infty}$ is an explicit very singular
solution. If $h(t)\geq e^{-\gw(t)/t}$ where $\gw$ is positive,
nondecreasing and satisfies
$$\int_{0}^1\frac{\sqrt{\omega(s)}ds}{s}=\infty,$$
then $u_{\infty}$ has a pointwise singularity at $(0,0)$. If the degeneracy of $h$ is stronger, namely
$$\liminf_{t\to 0}t\ln h(t)>-\infty,
$$
it is proved that the singularity of $u_{k}$ propagates along the axis $t=0$; at end, $u_{\infty}$ is nothing else than the (explicit) maximal solution $\Psi(t)$ of the O.D.E.
\begin{equation}\label{E1'''}
\Psi'+h( t)\Psi^p=0\quad\text { in }(0,\infty).
\end{equation}

A very general and probably difficult {\bf open problem} generalizing (\ref{E1}) and (\ref{E1''}) is to study the propagation phenomenon of singularities starting from $(0,0)$ when (\ref{E1}) is replaced by
\begin{equation}\label{E1^4}
\prt_tu-\Gd u+h( x,t)\abs u^{p-1}u=0\quad\text { in }Q_{T},
\end{equation}
where $h\in C(\overline Q_{T})$ is nonnegative and vanishes only on a curve $\Gg\subset \overline Q_{T}$ starting from $(0,0)$. It is expected that two types of phenomena should occur:\smallskip

\noindent (i) either $u_{\infty}$ has a pointwise singularity at $(0,0)$,\smallskip

\noindent (ii) or $u_{\infty}$ is singular along $\Gg$ or a connected part of $\Gg$ containing $(0,0)$.\smallskip

It is natural to conjecture that the order of degeneracy should be measured in terms of the parabolic distance to $\Gg$ and of the slope of $\Gg$ in the space $\BBR^N\ti\BBR$. This could serve as a starting model for nonlinear heat propagation in inhomogeneous fissured media.

\medskip

\noindent Our paper is organized as follows: 1 Introduction - 2
The power case -  3 Pointwise singularities - 4 Existence of razor
blades.\medskip

\noindent {\bf Acknowledgements } The authors have been supported
by INTAS grant Ref. No : 05-1000008-7921.
\section {The power case}
\setcounter{equation}{0}

In this section we assume that $h(x)=\abs x^{\gb}$ with $\gb\in\BBR$ and the equation under consideration
is the following
\begin{equation}\label{E2}
\prt_tu-\Gd u+\abs x^{\gb}\abs u^{p-1}u=0\quad\text { in }Q_{T}:=\BBR^N\ti (0,T)
\end{equation}
with $p>1$. By a solution we mean a function $u\in C^{2,1}(Q_{T})$. Let $E(x,t)=(4\gp t)^{-N/2}e^{-|x|^2/4t}$
be the heat kernel in $Q_{T}$ and $\BBE[\gf]$ the heat potential of a function (or measure) $\gf$ defined by
\begin{equation}\label{pc0}
\BBE[\gf](x,t)=\myfrac{1}{(4\gp t)^{N/2}}\myint{\BBR^N}{}
e^{-|x-y|^2/4t}\gf(y)\,dy.
\end{equation}
If there holds
\begin{equation}\label{pc1}
\dint_{Q_{T}}E^p(x,t)|x|^{\gb}dx\,dt<\infty,
\end{equation}
it is easy to prove (see \cite[Prop 1.2]{MV2}, and \cite[Th 6.12]{Ve}),  that for any $k\in\BBR$,
there exists a unique function $u=u_{k}\in L^1 (B_{R}\ti (0,T))\cap L^p(B_{R}\ti (0,T);|x|^\gb dx)$ such that
\begin{equation}\label{pc2}
\dint_{Q_{T}}\left(-u\prt_{t}\gz-u\Gd\gz+|x|^{\gb}\abs u^{p-1}u\gz\right)dx\,dt=k\gz(0,0),
\end{equation}
for any $\gz\in C^{2,1}_{0}(\BBR^N\ti [0,T))$. By the maximum principle
$k\mapsto u_{k}$ is increasing. Next, it is straightforward that (\ref{pc1}) is fulfilled as soon as
\begin{equation}\label{pc3}
\gb>\max\{N(p-1)-2;-N\}
\end{equation}
\subsection{The {\it a priori} estimate and the maximal solution}
In order to prove an {\it a priori} estimate, we introduce the auxiliary
$N$ dimensional equation in the variable $\eta=x/\sqrt t$
\begin{equation}\label{aux1}
-\Gd f-\myfrac{1}{2}\eta.\nabla f-\gamma f+\abs\eta^\gb\abs f^{p-1}f=0,
\end{equation}
where $\gamma=(2+\gb)/2(p-1)$.
\bprop{auxpr}Let $a>0$ and $\gb\in\BBR$; then there exists a unique nonnegative function $F_{a}\in H_{loc}^1(B_{a})\cap L_{loc}^{p+1}(B_{a};|\eta|^\gb d\eta)$ solution of (\ref{aux1}) and
 satisfying
\begin{equation}\label{aux2}
\lim_{|\eta|\to a}F_{a}(\eta)=\infty.
\end{equation}
Furthermore $a\mapsto F_{a}$ is decreasing.
\es
\Proof Set $K(\eta)=e^{\abs\eta^2/4}$. Then (\ref{aux1}) becomes
\begin{equation}\label{aux3}
-K^{-1}div (K\nabla f)-\gamma f+\abs\eta^\gb\abs f^{p-1}f=0.
\end{equation}
{\it Step 1- Boundary behaviour. } First we claim that
\begin{equation}\label{bh1}
\lim_{\abs \eta\to a}(a-\abs\eta)^{2/(p-1)}F_{a}(\eta)=\left(\myfrac{2(p+1)}{a^{p\gb}(p-1)^2}\right)^{1/(p-1)}.
\end{equation}
Actually, if $0<b<|\eta|<a$, $u$ satisfies
$$
-K^{-1}div (K\nabla F_{a})-\gamma F_{a}+ CF_{a}^p\leq 0
$$
with $C=\min\{a^\gb,b^\gb\}$. We perform a standard variant of
the two-sides estimate method used in \cite{Ve0} : we set
$\Gg:=B_{\gr}\setminus B_{b}$ with $b<\gr<a$, $\ga=(\gr-b)/2$ and
denote by $z$ the solution of
\begin{equation}\label{sol1}\left\{\BA {l}
z''-Cz^p=0\quad\text {in }(-\ga,\ga)\\
z(-\ga)=z(\ga)=\infty.
\EA\right.\end{equation}
Then $z$ is an even function and is computed by the formula
\begin{equation}\label{sol2}\myint{z(t)}{\infty}\myfrac{ds}{\sqrt{s^{p+1}-z(0)^{p+1}}}=
\sqrt{\myfrac{2C}{p+1}}(\ga-t)\forevery t\in [0,\ga).
\end{equation}
Notice also that $\lim_{\ga\to 0}z(t)=\infty$, uniformly on $(-\ga,\ga)$ and
\begin{equation}\label{sol0}
\lim_{t\to\ga}(t-\ga)^{2/(p-1)}z(t)=\left(\myfrac{2(p+1)}{C^p(p-1)^2}\right)^{1/(p-1)}.
\end{equation}
 We set $Z(\eta)=z(\abs\eta-(\gr+b)/2)$ and we look for a super-solution in $\Gg$ under the form $w=MZ(\eta)$ ($M>1$).
 Then
$$\BA {l}-K^{-1}div (K\nabla w)-\gamma w+ Cw^p
=M\left((M^{p-1}-1)Cz^p-\left(\myfrac{N-1}{\abs\eta}+\myfrac{\abs\eta}{2}\right)z'-\gamma z\right).
\EA$$
Since
$$z'(t)=\sqrt{\myfrac{2C}{p+1}}\sqrt{z^{p+1}(t)-z(0)^{p+1}}<C^*z^{(p+1)/2}(t),\quad\text {with }C^*=\sqrt{\myfrac{2C}{p+1}},
$$
we derive
\begin{equation}\label{sol3}-K^{-1}div (K\nabla w)-\gamma w+ Cw^p
\geq M\left((M^{p-1}-1)Cz^p-\left(\myfrac{N-1}{b}+\myfrac{a}{2}\right)C^*z^{(p+1)/2}-\gamma z\right)
\end{equation}
on $\{\eta:(\gr-b)/2\abs\eta<\gr\}$; and the same inequality holds true on
$\{\eta:\gr<\abs\eta<(\gr-b)/2\}$, up to interverting $a$ and $b$. For any $M>1$, we can choose $b>0$ such that for any $b<\gr<a$, the right-hand side of (\ref{sol3}) is positive and maximum principle applies in $B_{\gr}\setminus B_{b}$. Thus $MZ\geq F_{a}$ in $\Gg$. Furthermore, the previous comparison still holds if we take $\gr=a$, which implies $\ga=(a-b)/2$. Therefore, using the explicit value of $C$
\begin{equation}\label{bh1'}
\limsup_{\abs \eta\to a}(a-\abs\eta)^{2/(p-1)}F_{a}(\eta)\leq M\left(\myfrac{2(p+1)}{\min\{a^{p\gb},b^{p\gb}\}(p-1)^2}\right)^{1/(p-1)}.
\end{equation}
Because $M>1$ and $0<b<a$ are arbitrary, we derive
\begin{equation}\label{bh1''}
\limsup_{\abs \eta\to a}(a-\abs\eta)^{2/(p-1)}F_{a}(\eta)\leq \left(\myfrac{2(p+1)}{a^{p\gb}(p-1)^2}\right)^{1/(p-1)}.
\end{equation}
For the estimate from below we notice that $u$ satisfies
$$
-K^{-1}div (K\nabla F_{a})-\gamma F_{a}+ \tilde CF_{a}^p\geq 0
$$
in $\{\eta:b<\abs\eta<a\}$, with $\tilde C=\max\{a^\gb,b^\gb\}$. Taking now $\ga=a-b$, we denote by $\tilde z$ the positive solution of
\begin{equation}\label{sol6}\left\{\BA {l}
\tilde z''+\gamma \tilde z-\tilde C\tilde z^p=0\quad\text {in }(0,\ga)\\
\tilde z(0)=0\\
\tilde z(\ga)=\infty.
\EA\right.\end{equation}
Then $\tilde z$ is  computed by the formula
\begin{equation}\label{sol7}\myint{\tilde z(t)}{\infty}\myfrac{ds}
{\sqrt{\tilde z'\,\!\!^{2}(0)-\gamma s^2+2\tilde Cs^{p+1}/(p+1)}}=
\ga-t\forevery t\in [0,\ga),
\end{equation}
and formula (\ref{sol0}) is valid provided $C$ be replaced by $\tilde C$. We fix $A\in\prt B_{a}$ with coordinates $(a,0,...,0)$, and look for a subsolution under the form $\tilde w(\eta)=M\tilde z(\eta_{1}-b)$ with $0<M<1$. Then
$$
-K^{-1}div (K\nabla \tilde w)-\gamma \tilde w+ \tilde C\tilde w^p=
\tilde M\left((\tilde M^{p-1}-1)\tilde z^p-\myfrac{\eta_{1}}{2}\tilde w'\right)\leq 0,
$$
since $\tilde w'\geq 0$. Applying again the maximum principle, we derive $\tilde w(\eta)\leq F_{a}$ in $B_{a}\cap\{\eta:b<\eta_{1}<a\}$. But clearly the direction $\eta_{1}$ is arbitrary and can be replaced by any radial direction. Thus
\begin{equation}\label{bh10''}
\liminf_{|\eta|\to a}(a-|\eta|)^{2/(p-1)}F_{a}(\eta)\geq
\tilde M\left(\myfrac{2(p+1)}{\max\{a^{p\gb},b^{p\gb}\}(p-1)^2}\right)^{1/(p-1)}.
\end{equation}
In turn, (\ref{bh10''}) implies
\begin{equation}\label{bh11''}
\liminf_{\abs \eta\to a}(a-\abs\eta)^{2/(p-1)}F_{a}(\eta)\geq \left(\myfrac{2(p+1)}{a^{p\gb}(p-1)^2}\right)^{1/(p-1)},
\end{equation}
and (\ref{bh1}) follows from (\ref{bh1''}) and (\ref{bh11''}).
\medskip

\noindent{\it Step 2- Uniqueness. }  If $F'$ is another nonnegative solution of (\ref{aux1}) satisfying the same boundary blow-up conditions, then for any $\ge>0$, $F'_{\ge}=(1+\ge)F' $ is a super solution. Thus, for $\gd>0$,
$$\BA{l}
\dint_{B_{a}}\!\!\left(-\myfrac{div (K\nabla F_{a})}{F_{a}+\gd}
+\myfrac {div (K\nabla F'_{\ge})}{F'_{\ge}+\gd}+\abs\eta^\gb \left(\myfrac{F^{p}_{a}}{F_{a}+\gd}-\myfrac{F'_{\ge}\,\!^{p}}{F'_{\ge}+\gd}\right)K
\right)((F_{a}+\gd)^{2}-(F'_{\ge}+\gd)^{2})_{+}d\eta\\[4mm]
\phantom{-------}
\leq \gamma\dint_{B_{a}}\left(\myfrac{F_{a}}{F_{a}+\gd}-\myfrac{F'_{\ge}}{F'_{\ge}+\gd}\right)((F_{a}+\gd)^{2}-(F'_{\ge}+\gd)^{2})_{+}Kd\eta.
\EA$$
By monotonicity
$$\left(\myfrac{F^{p}_{a}}{F_{a}+\gd}-\myfrac{F'_{\ge}\,\!^{p}}{F'_{\ge}+\gd}
\right)((F_{a}+\gd)^{2}-(F'_{\ge}+\gd)^{2})_{+}\geq 0,
$$
and
$$0\leq\left(\myfrac{F_{a}}{F_{a}+\gd}-\myfrac{F'_{\ge}}{F'_{\ge}+\gd}\right)((F_{a}+\gd)^{2}-(F'_{\ge}+\gd)^{2})_{+}\leq ((F_{a}+\gd)^{2}-(F'_{\ge}+\gd)^{2})_{+}.
$$
By Lebesgue's theorem, since (\ref{bh1}) implies that $((F_{a}+\gd)^{2}-(F'_{\ge}+\gd)^{2})_{+}$ has compact support in $B_{a}$,
$$\lim_{\gd\to 0}\dint_{B_{a}}\left(\myfrac{F_{a}}{F_{a}+\gd}-\myfrac{F'_{\ge}}{F'_{\ge}+\gd}\right)((F_{a}+\gd)^{2}-(F'_{\ge}+\gd)^{2})_{+}Kd\eta=0.
$$
Using Green formula, we obtain
$$\BA {l}
\dint_{B_{a}}\!\!\left(-\myfrac{div (K\nabla F_{a})}{F_{a}+\gd}
+\myfrac {div (K\nabla F'_{\ge})}{F'_{\ge}+\gd}\right)((F_{a}+\gd)^{2}-(F'_{\ge}+\gd)^{2})_{+}Kd\eta\\[4mm]
\phantom{----}=
\dint_{F_{a}\geq F'_{\ge}}\left(\abs{\nabla F_{a}-\myfrac{F_{a}+\gd}{F'_{\ge}+\gd}\nabla F'_{\ge}}^2+\abs{\nabla F'_{\ge}-\myfrac{F'_{\ge}+\gd}{F_{a}+\gd}\nabla {F_{a}}}^2\right)Kd\eta\geq 0.
\EA$$
Letting $\gd\to 0$, we derive, by Fatou's theorem,
$$\dint_{F_{a}\geq F'_{\ge}}
\left(F^{p-1}_{a}-F'_{\ge}\,\!^{p-1}
\right)(F_{a}^{2}-F'_{\ge}\,\!^{2})Kd\eta\leq 0.
$$
Thus $F_{a}\leq F'_{\ge}$. Since $\ge$ is arbitrary, $F_{a}\leq F'$. The reverse inequality is the same. The monotonicity of $a\mapsto F_{a}$ is proved in a similar way, by the previous form of maximum principle.
\medskip

\noindent{\it Step 3- Existence with finite boundary value. } We
shall first prove the existence of a positive solution $w_{k}$ of
(\ref{aux1}) with boundary value equal to $k>0$  for small value
of $a$, and we shall let $k\to\infty$ in order to obtain one
solution satisfying (\ref{aux2}). We denote by $J_{a}$ the
functional defined over $H^1_{0}(B_{a})\cap
L^{p+1}(B_{a};|\eta|^\gb d\eta)$ by
$$ J_{a}(w)=\myfrac{1}{2}\myint{B_{a}}{}
\left(\abs{\nabla w}^2-\gamma w^2+\myfrac{1}{p+1}|\eta|^{\gb}|w|^{p+1}\right)K(\eta)d\eta.
$$
Let $k>0$ and $\gk\in C^1(\overline B_{a})$ with $0\leq \gk(\eta)\leq k$, $supp (\gk)\subset \overline B_{a}\setminus B_{a/2}$, $\gk(\eta)\equiv k$ on $\overline B_{a}\setminus B_{2a/3}$. If
$v\in H^1_{0}(B_{a})\cap L^{p+1}(B_{a};|\eta|^\gb d\eta)$ and
$w:=v+\gk$, then
$$J_{a}(w)=J_{a}(v+\gk)\geq J_{a}(v)+J_{a}(\gk)
+\myint{B_{a}}{}\left(\nabla v.\nabla\gk-\gamma v\gk-|\eta|^\gb|v|^p\gk\right)K(\eta)d\eta.
$$
Since $\gamma\leq \gl_{a}$, it follows from Cauchy-Schwarz and H\"older-Young inequalities that
$$J_{a}(w)\geq (1-\ge^2)J_{a}(v)-\myfrac{p^p}{\ge^{2p}}J_{a}(\gk)
$$
for $0<\ge<1$. Because $\lim_{a\to 0}\gl_{a}=\infty$, there exists $a_{0}\in (0,\infty]$ such that, for any $0<a<a_{0}$, $J_{a}(v)$ is bounded from below on
$H^1_{0}(B_{a})\cap L^{p+1}(B_{a};|\eta|^\gb d\eta)$. Thus there exists a minimizer $w_{k}$ such that $w_{k}=v+\gk$ with $v$ in the above space;  $w_{k}$ is a solution of (\ref{aux1}) and $w_{k}|_{\prt B_{a}}=k$. Furthermore $w_{k}$ is positive. Notice that if $\gamma\leq 0$, $a_{0}=\infty$, in which case there exists a solution $w_{k}$ for any $k>0$ and any $a>0$. The uniqueness of $w_{k}>0$, is a consequence of the monotonicity of the mapping $k\mapsto w_{k}$ that we prove by a similar argument as in Step 2: if  $k<k'$, there holds
$$\BA {l}
\dint_{w_{k}>w_{k'}}\left(w_{k}^{p-1}-w_{k'}^{p-1})(w^2_{k}-w_{k'}^{2})
\right)\abs\eta^\gb K d\eta
\leq 0,
\EA
$$
which implies $w_{k}<\tilde w_{k}$. Uniqueness and radiality follows immediately, thus $w_{k}$ solves the differential equation
\begin{equation}\label{ode}\left\{\BA {l}
-w''-\left(\myfrac{N-1}{r}+\myfrac{r}{2}\right)w'-\gamma w+r^\gb w^p=0\quad\text{on }(0,a)\\
w(a)=k\quad\mbox{and }w\in H_{rad}^1(B_{a})\cap L^{p+1}_{rad}(B_{a};|\eta|^\gb d\eta).
\EA\right.
\end{equation}
Next we shall assume $\gamma >0$, equivalently $\gb>-2$. If
$w_{k}$ is a positive solution of (\ref{ode}) and $\gl>1$ (resp.
$\gl<1$) $\gl w_{k}$ is a super-solution (resp. a sub-solution)
larger (resp. smaller) than $w_{k}$. Note that $\gb>-2$ implies $w_{k}(0)>0$ while $\gb>-1$ implies also $w'_{k}(0)=0$. Thus, by \cite {RRV}, there
exists a solution $w_{\lambda k}$ with boundary data $\lambda k$,
and this solution is positive because $w_{k}\leq w_{\lambda k}\leq
\lambda w_{ k}$ (resp. $\lambda w_{k}\leq w_{\lambda k}\leq  w_{
k}$). Consequently, the set $\CA$ of the positive $\tilde a$ such
that there exists a positive solution of (\ref{ode}) on $(0,a)$
for any $a<\tilde a$ is not empty and independent of $k$.
Furthermore, if for some $\tilde a>0$ and some $k_{0}>0$,  there
exists some positive $w_{k_{0}}$ solution of (\ref{ode}) on
$0,\tilde a)$, then for any $0<a<\tilde a$ and any $k>0$, there
exists a positive solution $w_{k}$ of (\ref{ode}). Since $r\mapsto
\max\{k,(\gamma_{+}a^{-\gb})^{1/(p-1)}\}$ is a super-solution,
there holds
\begin{equation}\label{ode1}
w_{k}(r)\leq \max\{k,(\gamma_{+}a^{-\gb})^{1/(p-1)}\}\forevery r\in [0,a].
\end{equation}
Let us assume that $a^*=\sup \CA<\infty$. Because of (\ref{ode1}) and local regularity of solutions of elliptic equations, for any $\ge,\ge'>0$, $w'_{k}(a)$ is bounded uniformly with respect if $ \ge\leq a<a^*-\ge'$. But since (\ref{ode}) implies
$$a^{N-1}e^{a^2/4}w'_{k}(a)=\ge^{N-1}e^{\ge^2/4}
w'_{k}(\ge)+\myint{\ge}{a}(r^\gb w^p_{k}-\gamma w_{k})r^{N-1}e^{r^2/4}dr,
$$
$w'_{k}(a)$ is actually uniformly bounded on $ [\ge,a^*)$.
It follows from the local existence and uniqueness theorem that there exists $\gd>0$, independent of $a<a^*$ such that there exists a unique solution $z$ defined on $[a,a+\gd]$ to
\begin{equation}\label{ode2}\left\{\BA {l}
-z''-\left(\myfrac{N-1}{r}+\myfrac{r}{2}\right)z'-\gamma z+r^\gb z^p=0\quad\text{on }(0,a)\\
z(a)=k,\,z'(a)=w'_{k}(a),
\EA\right.
\end{equation}
and $\gd$ and $k>0$ can be chosen such that $z>0$ in $[a,a+\gd]$.
This leads to the existence of a positive solution to (\ref{ode})
on $[0,a+\gd]$. If $a^*-a<\gd$, which contradicts the maximality
of $a^*$. Therefore $a^*=\infty$.\medskip

\noindent{\it Step 4- End of the proof. }We have already seen that $k\mapsto w_{k}$ is increasing. By Step 1, we know that, for any $a>0$, and some $b<a$, there holds
\begin{equation}\label{ode'1}
w_{k}(|\eta|)\leq C(a-\abs\eta)^{-1/(p-1)}\text {on }B_{a}\setminus B_{b}.
\end{equation}
In particular
$$w_{k}(b)\leq C^*=C^*(a,b,p,N)
$$
Next
\begin{equation}\label{ode''1}
w_{k}(r)\leq \max\{C^*,(\gamma_{+}b^{-\gb})^{1/(p-1)}\}\forevery r\in [0,b].
\end{equation}
Combining (\ref{ode'1}) and (\ref{ode''1}) implies that $w_{k}$ is locally uniformly bounded on $[0,a)$. Since $k\mapsto w_{k}$ is increasing, the existence of $F_{a}:=w_{\infty}=\lim_{k\to\infty}w_{k}$ follows. The fact that $a\mapsto F_{a}$ decreases is a consequence of the fact that $F_{a'}$ is finite on $\prt B_{a}$ for any $a<a'$.
\qeda
\medskip

\noindent\Remark
In the sequel we set $F_{\infty}=\lim_{a\to\infty}F_{a}$. Then $F_{\infty}$ is a nondecreasing, nonnegative solution of (\ref{aux1}).
Using asymptotic analysis, is is easy to prove that  there holds:\smallskip

\noindent (i) if $\gb\neq 0$
\begin{equation}\label{asym}
F_{\infty}(\eta)=\left(\myfrac{1}{p-1}\right)^{1/(p-1)}|\eta|^{-\gb/(p-1)}(1+\circ(1))\quad\text{as }|\eta|\to\infty;
\end{equation}

\noindent (ii) if $\gb=0$,
\begin{equation}\label{pc5^}
F_{\infty}(\eta)\equiv\left(\myfrac{1}{p-1}\right)^{1/(p-1)}.
\end{equation}

Furthermore, if $\gb>-2$, it follows by the strict maximum principle that $F_{a}(0)=\min\{F_{a}(\eta):\abs\eta<a\}>0$. This observation plays a fundamental role for obtaining estimate from above.
\bprop {KOpr} Assume $p>1$ and $\gb>-2$. Then any solution $u$ of (\ref{E2}) in $Q_{T}$ which verifies
\begin{equation}\label{pc4}
\lim_{t\to 0}u(x,t)=0\quad\forall x\neq 0,
\end{equation}
satisfies
\begin{equation}\label{pc5}
|u(x,t)|\leq \min\left\{c^*\abs x^{-(2+\gb)/(p-1)};t^{-(2+\gb)/2(p-1)}F_{\infty}(x/\sqrt t)\right\}\quad\forall (x,t)\in Q_{T}\setminus\{0\},
\end{equation}
where $c^*=c^*(N,p,\gb)$.\es
\Proof Let $\ge>0$ and $a>0$ and
$\CP_{a,\ge}=\{(x,t):t>\ge,|x|/\sqrt{t-\ge}<a\}$. By the previous remark $\min F_{a}>0$, thus the function
$W(x,t)=(t-\ge)^{-(2+\gb)/2(p-1)}F_{a}(|x|/\sqrt{t-\ge})$, which is a solution of (\ref{E2}) in $\CP_{a,\ge}$ tends to infinity on the boundary on $\CP_{a,\ge}$;  since $u$ is finite in $Q_{T}\cap \CP_{a,\ge}$, $W$ dominates $u$ in this domain. Letting successively $\ge\to 0$ and $a\to\infty$ yields to $u\leq F_{\infty}$. The estimate from below is similar. Next we consider $x\in \BBR^N\setminus\{0\}$, then $v=|u|$ satisfies (by Kato's inequality)
$$\prt_{t}v-\Gd v+C(x)v^p\leq 0\quad\text{in }B_{|x|/2}(x)\ti (0,T),
$$
where $C(x)=\max\{(|x|/2)^\gb
;(3|x|/2)^\gb\}$. It is easy to construct a function under the form
$w(y)=\Gl\left(|x|^2-4|x-y|^2\right)^{-2/(p-1)}$ which satisfies
$$\left\{\BA {l}-\Gd w+C(x)w^p=0\quad\text {in }B_{|x|/2}(x)\\
\lim_{|x-y|\to|x|/2}w=\infty,
\EA\right.$$
with $\Gl=\Gl(x)=c^*\abs x^{(2-\gb)/(p-1)}$, $c^*=c^*(N,p,\gb)>0$. Using (\ref{pc4}), it follows from Lebesgue's theorem that $u(y,t)\leq w(y)$ in $B_{|x|/2}(x)\ti [0,T)$, thus $u(x,t)\leq w(x)=c^*|x|^{-(2+\gb)/(p-1)}$. Estimate from below is similar.
\qeda\medskip

The construction of the first part of the proof of \rprop{KOpr} (estimate in
$\CP_{a,\ge}$) shows that, without condition (\ref{pc4}), equation  (\ref{E2}) admits a maximal solution $u_{M}$.
\bprop{maxsol}Assume $p>1$ and $\gb>-2$.
Then any solution $u$ to (\ref{E2}) satisfies
\begin{equation}\label{pc5*}
|u(x,t)|\leq u_{M}(x,t):= t^{-(2+\gb)/2(p-1)}F_{\infty}(x/\sqrt t)\quad\forall (x,t)\in Q_{T}\setminus\{0\}.
\end{equation}
\es

As a variant of (\ref{pc5}), we have the following Keller-Osserman type parabolic estimate which extends the classical one due to Brezis and Friedman in the case $\gb=0$ (see \cite{BF}).
\bprop{KOpr'}Under the assumptions of \rprop {KOpr} there holds
\begin{equation}\label {pc5'}
|u(x,t)|\leq \myfrac{\tilde c}{(\abs x^2+t)^{(2+\gb)/2(p-1)}}\quad\forall (x,t)\in Q_{T}\setminus\{0\},
\end{equation}
with $\tilde c=\tilde c(N,p,\gb)$.
\es
\Proof Assume $\abs x^2\leq t$, then
\begin{equation}\BA {l}\label{pc5''}
\myfrac{1}{(\abs x^2+t)^{(2+\gb)/2(p-1)}}\geq 2^{-(2+\gb)/2(p-1)}t^{-(2+\gb)/2(p-1)}\\[2mm]\phantom{\myfrac{1}{(\abs x^2+t)^{(2+\gb)/2(p-1)}}}
\geq \myfrac{2^{-(2+\gb)/2(p-1)}}{\min\{F_{\infty}(\eta):\abs\eta\leq 1\}}t^{-(2+\gb)/2(p-1)}F_{\infty}(x/\sqrt t).
\EA\end{equation}
Assume $\abs x^2\geq t$, then
\begin{equation}\BA {l}\label{pc5'''}
\myfrac{1}{(\abs x^2+t)^{(2+\gb)/2(p-1)}}\geq 2^{-(2+\gb)/(p-1)}|x|^{-(2+\gb)/(p-1)}.
\EA\end{equation}
Combining (\ref{pc5''}) and (\ref{pc5'''}) gives (\ref{pc5'}).\qeda
\medskip

\subsection{Isolated singularities and the very singular solution}
\bth{remov} Assume $p>1$ and $-2<\gb\leq N(p-1)-2$. Then any
solution $u$ to (\ref{E2}) which satisfies (\ref{pc4}) is
identically $0$. \es
\Proof  If $-(2+\gb)/(p-1)+N-1>-1$,
equivalently $\gb<N(p-1)-2$, the function $x\mapsto \abs
x^{-(2+\gb)/(p-1)}$ is locally integrable in $\BBR^N$, thus
$u(.,t)\to 0$ in $L^1_{loc}(\BBR^N)$ as $t\to 0$. For $\ge>0$
there exists $R=R(\ge)$ such that $u(x,t)\leq \ge$ for any $\abs
x\geq R$ and $t>0$. Thus
\begin{equation}\label{pc6}
u(x,t+\gt)\leq\ge+\BBE[u\chi_{_{B_{R}}}u(.,\gt)](x,t)\forevery t>0,\gt>0\,\text{and }x\in\BBR^N,
\end{equation}
where $\BBE[\gf]$ denotes the heat potential of the measure $\gf$ (see (\ref{pc0})). Letting successively $\gt\to 0$ and $\ge\to 0$, yields to $u\leq 0$. In the same way $u\geq 0$. In the case $\gb=N(p-1)-2$ estimate (\ref{pc5'}) reads
$$|u(x,t)|\leq \myfrac{\tilde c}{(\abs x^2+t)^{N/2}}.
$$
From this estimate, the proof of \cite[Th 2, Steps 5, 6]{BF} applies and we recall briefly the steps\smallskip

\noindentÊ(i) By choosing positive test functions $\gf_{n}$ which vanish in $\CV_{n}=\{(x,t):\abs x^2+t\leq n^{-1}\}$ and are constant on $\CV'_{n}=\{(x,t):\abs x^2+t\geq 2n^{-1}\}$, we first prove that, for any $\gr>0$,
\begin{equation}\label{pc7}
\dint_{B_{\gr}\ti (0,T)}\left(|u(x,t)|+|x|^{\gb}|u|^p\right)dxdt<\infty.
\end{equation}
Thus, using the same test function, we derive that the identity
\begin{equation}\label{pc8}
\dint_{Q_{T}}\left(-u\prt_{t}\gz-u\Gd\gz+|x|^{\gb}\abs u^{p-1}u\gz\right)dx\,dt=0,
\end{equation}
holds for any $\gz\in C^{2,1}_{0}(\BBR^N\ti [0,T))$. The uniqueness yields to
$u=0$.
\qeda
\medskip

\noindent{\it Proof of Theorem A- case I. } In the case $-2<\gb\leq N(p-1)-2$, the result is a consequence of \rth{remov}. Next we assume $\gb\leq -2$. If $f$ is a solution of (\ref{int4}), it satisfies
$$f(\eta)=\circ(|\eta|^{-(2+\gb)/(p-1)})\quad\text {as }|\eta|\to \infty.
$$
If $\gb=-2$, the equation becomes
$$-\Gd f-\myfrac{1}{2}\eta.\nabla f+|\eta|^{-2}|f|^{p-1}f=0,
$$
and $f(\eta)\to 0$ at infinity. Since any positive constant is a supersolution, $f\leq 0$. Similarly $f\geq 0$.\smallskip

\noindent If $\gb<-2$, for $\ge>0$ the function $\eta\mapsto \ge|\eta|^{-(2+\gb)/(p-1)}=\psi (\eta)$ belongs to $ W^{1,1}_{loc}(\BBR^N)$ since $\gb<-2$ and satisfies
$$\BA {l}
-\Gd \psi-\myfrac{1}{2}\eta.\nabla \psi-\myfrac{2+\gb}{2(p-1)}\psi+|\eta|^{\gb}|\psi|^{p-1}\psi\\\phantom{-----------}
=\ge r^{-(2+\gb)/(p-1)-2}\left(\left(\myfrac{2+\gb}{p-1}\right)\left(\myfrac{2+\gb}{p-1}+2-N\right)+\ge^{p-1}\right).
\EA$$
 Therefore, either if $N\geq 2$ or $N=1$ and $\gb\leq -(p+1)$, $\psi$ is a super-solution of (\ref{int4}) for any $\ge>0$. The conclusion follows as above.\smallskip

 \noindent Finally we treat the case $N=1$ and $-(p+1)<\gb<-2$ where there exists a particular solution of
 $$f''+\myfrac{r}{2}f'+\myfrac{2+\gb}{2(p-1)}f-r^{\gb}|f|^{p-1}f=0\quad\mbox {on }\BBR_{+},
 $$
 under the form $f_{1}(r)=A_{\gb,p}r^{-(2+\gb)/(p-1)}$. Furthermore, if $f\geq 0$ (which can be always assumed by the maximum principle), it is a subsolution of the linear equation
 $$\gf''+\myfrac{r}{2}\gf'+\myfrac{2+\gb}{2(p-1)}\gf=0
 $$
Noticing that this equation has a solution $\gf_{1}$ which has the same  behaviour at infinity than the explicit solution of (\ref{int2}), namely
$$\gf_{1}(r)=cr^{-(2+\gb)/(p-1)}(1+\circ(1)),
$$
by standard methods (see e.g. \cite [Prop A1]{MV0}), the second solution $\gf_{2}$ behaves in the following way
$$\gf_{2}(r)=cr^{(2+\gb)/(p-1)-1}e^{-r^2/4}(1+\circ(1))\quad\mbox{as }r\to\infty.
$$
Consequently, by the maximum principle, any solution $f$ of (\ref{int2}) on $\BBR$ such that $f(r)=\circ(\gf_{1}(r))$ at infinity, verifies
\begin{equation}\label {R1}
|f(r)|\leq C|r|^{(2+\gb)/(p-1)-1}e^{-r^2/4}\quad\mbox{for }|r|\geq 1.
\end{equation}
Using the equation, we obtain that
$$f'(r)=e^{r^2/4}\myint{r}{\infty}\left(s^{\gb}|f(s)|^{p-1}f(s)-\myfrac{2+\gb}{p-1}f(s)\right)ds,
$$
thus
\begin{equation}\label {R2}|f'(r)|\leq Cr^{(2+\gb)/(p-1)-2}e^{-r^2/4}\quad\mbox{for }|r|\geq 1.
\end{equation}
Since $f\in H^1_{loc}(\BBR)$, we derive that for any $n\in\BBN_{*}$,
$$\myint{-n}{n}\left(f'\,\!^{2}-\myfrac{2+\gb}{p-1}f^2\right)e^{r^2/4}dr
\leq e^{n^2/4}\left(f(n)f'(n)-f(-n)f'(-n)\right).$$
Because of (\ref{R1}) and (\ref{R1}), this last term tends to $0$ as $n\to\infty$. Therefore
$$\myint{-\infty}{\infty}\left(f'\,\!^{2}-\myfrac{2+\gb}{p-1}f^2\right)e^{r^2/4}dr=0\Longrightarrow f=0,
$$
which end the proof.\qeda\medskip

\noindent\Remark The method of proof used in the case $N=1$ and $-p-1<\gb<-2$ is actually valid in any dimension, for any $\gb\leq-2$. But it relies strongly on the fact that $f\in H^1_{loc}$, while the other methods use only $f\in W^{1,1}_{loc}(\BBR^N).$

\bprop {u-infin}Assume $\gb> \max\{N(p-1)-2;-N\}$ . Then for any
$k>0$ there exists a unique solution $u_{k}$ of (\ref{E2}) with initial data $k\gd_{0}$. Furthermore $k\mapsto u_{k}$ is increasing and $u_{\infty}:=\lim_{k\to\infty}u_{k}$ satisfies $u_{\infty}(x,t)=t^{-(2+\gb)/2(p-1)}f_{\infty}(x/\sqrt t)$, where $f_{\infty}$ is positive, radially symmetric and satisfies
\begin{equation}\label{self}\left\{\BA {l}
-\Gd f_{\infty}-\myfrac{1}{2}\eta.\nabla f_{\infty}-\gamma f_{\infty}
+\abs\eta^{\gb}f_{\infty}^p=0\quad\text {in }\BBR^N\\\phantom{---}
\lim_{\abs\eta\to\infty}|\eta|^{(2+\gb)/(p-1)}f_{\infty}(\eta)=0.
\EA\right.\end{equation}
\es
\Proof The existence of  $u_{k}$ and the monotonicity of $k\mapsto u_{k}$ has already been seen. By the uniform continuity of the $u_{k}$ in any compact subset of $\bar Q_{T}\setminus\{(0,0)\}$, the function $u_{\infty}$ satisfies
\begin{equation}\label{self0}\lim_{t\to 0}u_{\infty}(x,t)=0\quad\forall x\neq 0.
\end{equation}
For $\ell>0$ and $u$ is defined in $Q_{\infty}$, we set
\begin{equation}\label{self1}T_{\ell}[u](x,t):=\ell^{(2+\gb)/2(p-1)}u(\sqrt \ell x,\ell t).
\end{equation}
If $u$ satisfies equation (\ref{E2}) in $Q_{\infty}$, $T_{\ell}[u]$ satisfies it  too. Because of uniqueness
\begin{equation}\label{self2}T_{\ell}[u_{k}]=u_{\ell^{(2+\gb)/2(p-1)-N/2}k}.
\end{equation}
Using the continuity of $u\mapsto T_{\ell}[u]$ and the definition of $u_{\infty}$, we can let $k\to\infty$ in (\ref{self2}) and derive (by taking $\ell t=1$ and replacing $t$ by $\ell$),
\begin{equation}\label{self3}T_{\ell}[u_{\infty}]=u_{\infty}
\Longrightarrow u_{\infty}(x,t)=t^{-(2+\gb)/2(p-1)}u_{\infty}(x/\sqrt t,1).
\end{equation}
Setting $f_{\infty}(\eta)=u_{\infty}(x/\sqrt t,1)$ with $\eta=x/\sqrt t$, it is straightforward that $f_{\infty}$ satisfies (\ref{self}) (using in particular \ref{self0}). Furthermore $f_{\infty}$ is radial and positive as the $u_{k}$ are.\qeda

\blemma{asymp} The function $f_{\infty}$ satisfies
\begin{equation}\label{asym1}
f_{\infty}(\eta)=c|\eta|^{2\gamma-N}e^{-|\eta|^2/4}\left(1+\circ(|\eta|^{-2})\right)\quad\text
{as }|\eta|\to\infty,
\end {equation}
for some $c=c_{N,p,\gb}>0$. Furthermore
\begin{equation}\label{asym1'}
f'_{\infty}(\eta)=-\myfrac{c}{2}c|\eta|^{2\gamma+1-N}e^{-|\eta|^2/4}\left(1+\circ(|\eta|^{-2})\right)\quad\text
{as }|\eta|\to\infty.
\end {equation}\es
\Proof Set $r=\abs\eta$ and denote $f_{\infty}(\eta)=f_{\infty}(r)$. Then $f_{\infty}$ satisfies,
\begin{equation}\label{asym2}
f_{\infty}''+\left(\myfrac{N-1}{r}+\myfrac{r}{2}\right)f_{\infty}'
+\gamma f_{\infty}-r^{\gb}|f_{\infty}|^{p-1}f_{\infty}=0\quad\text {on }(0,\infty),
\end {equation}
and $\lim_{r\to\infty}r^{2\gamma}f_{\infty}(r)=0$. We consider the auxiliary equation
\begin{equation}\label{auxe}
f''+\left(\myfrac{N-1}{r}+\myfrac{r}{2}\right)f'
+\gamma f=0\quad\text {on }(0,\infty).
\end {equation}
By \cite[Prop A1]{MV0}, (\ref{auxe}) admits two linearly independent solutions defined on $(0,\infty)$, $y_{1}$ and $y_{2}$ such that
\begin{equation}\label{auxe1}
y_{1}(r) = r^{-2\gamma} (1 + \circ(1))\quad\text{ and }\; y_{2}(r) = r^{2\gamma-N}e^{-r^2/4}(1 + \circ(1)),
\end {equation}
as $r\to\infty$. Next we choose $R>0$ large enough so that the maximaum principle applies for equation (\ref{auxe}) on $[R,\infty)$ and the $y_{j}$ are positive on the same interval. For $\gd>0$, $Y_{\gd}=\gd y_{1}+f_{\infty}(R)y_{2}/y_{2}(R)$ is a supersolution for (\ref{asym2}). Furthermore $f_{\infty}(r)=\circ (Y_{\gd})$ at infinity. Letting $\gd\to 0$ yields to
\begin{equation}\label{auxe2}
f_{\infty}(r)\leq \myfrac{f_{\infty}(R)}{y_{2}(R)}y_{2}(r)\forevery r\geq R.
\end {equation}
Using (\ref{auxe1}) we derive
$$0\leq f_{\infty}(\eta)\leq C|\eta|^{2\gamma-N}e^{-|\eta|^2/4}\forevery \abs\eta\geq 1.
$$
Plugging this estimate into (\ref{asym2}), we derive (\ref{asym1}) from standard perturbation theory for second order linear differential equation \cite[p. 132-133]{Bel}. Finally, (\ref{asym1'}) follows directly from (\ref{asym1}) and (\ref{asym2}).\qeda
\medskip

An alternative proof of the existence of $f_{\infty}$ is linked to calculus of variations. In the case $\gb=0$, this was performed by Escobedo and Kavian \cite{EK}. This construction is based upon the study of the following functional
\begin{equation}\label{J}
J(v)=\myfrac{1}{2}\myint{\BBR^N}{}\left(\abs{\nabla v}^2-
\gamma v^2+\myfrac{2}{p+1}\abs \eta^{\gb}\abs v^{p+1}\right)K(\eta)d\eta,
\end {equation}
defined over the functions in $H^1_{K}(\BBR^N)\cap L_{\abs\eta^\gb K}^{p+1}(\BBR^N)$.
\bprop{eke} Assume $p>1$ and $\gb> N(p-1)-2$. Then there exists a positive function
$\tilde f_{\infty}\in H^1_{K}(\BBR^N)\cap L_{\abs\eta^\gb K}^{p+1}(\BBR^N)$ satisfying
\begin{equation}\label{J'}
-\Gd \tilde f_{\infty}-\myfrac{1}{2}\eta.\nabla\tilde
f_{\infty}-\gamma\tilde f_{\infty} +|\eta|^\beta\tilde
f_{\infty}^p=0\quad\text{in }\BBR^N.
\end {equation}
\es

We recall that the eigenvalues of $-K^{-1}div(K\nabla.)$ are the $\gl_{k}=(N+k)/2$, with $k\in\BBN$ and the eigenspaces $H_{k}$ are generated by $D^\ga\gf$ where $\gf(\eta)=K^{-1}(\eta)=e^{-\abs\eta^2/4}$ and $\abs\ga =k$. It is straightforward to check that $J$ is $C^1$. In order to apply Ekeland Lemma, we have just to prove that $J$ is bounded from below in $H^1_{K}(\BBR^N)$. As we shall see it later on, the proof is easy when  $\gb<N(p-1)/2$, and more difficult when $\gb\geq N(p-1)/2$.
\blemma{est}For any $v\in H^1_{K}(\BBR^N)$, there holds
$$\myfrac{1}{4}\myint{\BBR^N}{}\left(2N+|\eta|^2\right)v^2K(\eta)d\eta\leq \myint{\BBR^N}{}|\nabla v|^2K(\eta)d\eta.
$$
\es
\Proof We borrow the proof to Escobedo and Kavian. Put $w=v\sqrt {K}$. Then
$$\sqrt {K}\nabla v=\nabla w-\myfrac{w}{2}\eta.
$$
Hence
$$\myint{\BBR^N}{}|\nabla v|^2K(\eta)d\eta=
\myint{\BBR^N}{}\left(|\nabla w|^2-w\nabla w.\eta+\myfrac{1}{4}w^2\abs\eta^2\right)d\eta.
$$
Because
$$-\myint{\BBR^N}{}w\nabla w.\eta d\eta=\myfrac{N}{2}
\myint{\BBR^N}{}w^2 d\eta,
$$
there holds
$$\myint{\BBR^N}{}|\nabla v|^2K(\eta)d\eta=
\myint{\BBR^N}{}\left(|\nabla w|^2+\myfrac{N}{2}w^2+\myfrac{1}{4}w^2\abs\eta^2\right)d\eta.
$$
This implies the formula.\qeda
\blemma{est2}Let $p>1$ and $\gb<N(p-1)/2$. For any $\ge>0$ there exists
$C=C(\ge,p)>0$ and $R=R(\ge,p)>0$ such that
$$\myint{\BBR^N}{}v^2K(\eta)d\eta\leq
\ge\myint{\BBR^N}{}|\nabla v|^2K(\eta)d\eta+
C\left(\myint{\BBR^N}{}|v|^{p+1}|\eta|^\gb K(\eta)\right)^{2/p+1}.
$$
\es
\Proof For $R>0$ there holds
$$\myint{|\eta|\leq R}{}v^2K(\eta)d\eta
\leq \left(\myint{|\eta|\leq R}{}|v|^{p+1}|\eta|^\gb K(\eta)d\eta\right)^{2/(p+1)}\left(\myint{|\eta|\leq R}{}|\eta|^{-2\gb/(p-1)}K(\eta)d\eta\right)^{(p-1)/(p+1)}.
$$
Since $\gb<N(p-1)/2\Longleftrightarrow N>2\gb/(p-1)$, we obtain
$$\left(\myint{|\eta|\leq R}{}|\eta|^{-2\gb/(p-1)}K(\eta)d\eta\right)^{(p-1)/(p+1)}=C(R,N,p).
$$
By \rlemma {est}
$$\myint{|\eta|\geq R}{}v^2K(\eta)d\eta
\leq \myfrac{4}{R^2}\myint{\BBR^N}{}|\nabla v|^2K(\eta)d\eta.
$$
The estimate follows by taking $\ge=4R^{-2}$.\qeda\medskip

It follows from the previous Lemmas that $J$ is bounded from below in the space $H^1_{K}(\BBR^N)\cap L_{\abs\eta^\gb K}^{p+1}(\BBR^N)$ whenever $N(p-1)/2-2<\gb<N(p-1)/2$. Next we consider the case $\gb>0$ and we shall restrict the study to radial functions.

\blemma{rad-decrease} Assume $\gb>0$. The functional $J$ is
bounded from below on the set
$$X=\left\{v\in H^1_{K}(\BBR^N)\cap L_{\abs\eta^\gb K}^{p+1}(\BBR^N): v\geq 0, v\text { radial and decreasing } \right\}.
$$
\es
\Proof For $0<\gd<R$, we write $J(v)=J_{\gd,R}(v)+J'_{\gd,R}(v)+J''_{\gd,R}(v)$ where
$$J_{\gd,R}(v)=\myfrac{1}{2}\myint{|\eta|\leq \gd}{}\left(\abs{\nabla v}^2
-\gamma v^2+\myfrac{2}{p+1}\abs \eta^{\gb}\abs v^{p+1}\right)K(\eta)d\eta,
$$
$$J'_{\gd,R}(v)=\myfrac{1}{2}\myint{\gd<|\eta|<R}{}\left(\abs{\nabla v}^2
-\gamma v^2+\myfrac{2}{p+1}\abs \eta^{\gb}\abs v^{p+1}\right)K(\eta)d\eta,
$$
and
$$J''_{\gd,R}(v)=\myfrac{1}{2}\myint{|\eta|>R}{}\left(\abs{\nabla v}^2
-\gamma v^2+\myfrac{2}{p+1}\abs \eta^{\gb}\abs v^{p+1}\right)K(\eta)d\eta.
$$

Using \rlemma{est2}, we fix $R$ large enough so that $J''_{\gd,R}$ is bounded from below in $H^1_{K}(\BBR^N)\cap L_{\abs\eta^\gb K}^{p+1}(\BBR^N)$.
By H\"older's inequality $J'_{\gd,R}$ is bounded from below, thus we are left with $J_{\gd,R}$. We assume that $v$ is positive, radial, nonincreasing and $v(\gd)=c=\min\{v(x):|x|\leq \gd\}$. Then
$$\abs v^{p+1}=v^{p+1}=(v-c+c)^{p+1}\geq (v-c)^{p+1}+c^{p+1}\;\text { and }\;v^2\leq 2(v-c)^2+2c^2,
$$
$$J_{\gd,R}\geq\myfrac{1}{2}\myint{|\eta|\leq \gd}{}\left(\abs{\nabla (v-c)}^2
-2\gamma (v-c)^2+\myfrac{2}{p+1}\abs \eta^{\gb}\abs {v-c}^{p+1}\right)K(\eta)d\eta +L(c),
$$
where
$$L(c)=\myfrac{c^{p+1}}{p+1}\myint{|\eta|\leq \gd}{}\abs \eta^{\gb}K(\eta)d\eta-\gamma c^2\myint{|\eta|\leq \gd}{}K(\eta)d\eta.
$$
Clearly $L(c)\geq M$ for some $M$ independent of $c$. Therefore we are reduced to study the functional $J_{\gd,R}$ defined by
$$J_{\gd,R}(w)=\myfrac{1}{2}\myint{|\eta|\leq \gd}{}\left(\abs{\nabla w}^2
-2\gamma w^2+\myfrac{2}{p+1}\abs \eta^{\gb}\abs {w}^{p+1}\right)K(\eta)d\eta
$$
over $H^1_{0, K}(B_{\gd})\cap L_{\abs\eta^\gb K}^{p+1}(B_{\gd})$. Here we can fix $\gd>0$ small enough so that the first eigenvalue of
$-K^{-1}div (K\nabla .)$ is larger than $2\gamma$, thus $J_{\gd,R}(v)$ is bounded from below in the class of radially symmetric nonincreasing, nonnegative functions $v$, and so is $J$.\qeda
\blemma {min}Let $v$ be a radially symmetric function in $H^1_{K}(\BBR^N)\cap L_{\abs\eta^\gb K}^{p+1}(\BBR^N)$. Then there exists a radially symmetric decreasing function $\tilde v\in H^1_{K}(\BBR^N)\cap L_{\abs\eta^\gb K}^{p+1}(\BBR^N)$ such that $J(\tilde v)\leq J(v)$.
\es
\Proof We define the two curves
$$C_{1}=\left\{(s,x)\in\BBR_{+}\ti\BBR_{+}:-2^{-1}\gamma x^2+(p+1)^{-1}s^{\gb}x^{p+1}=0\right\}=\left\{x=\left(2^{-1}(p+1)\gamma s^{-\gb}\right)^{1/(p-1)}\right\},
$$
and
$$C_{2}=\left\{(s,x)\in\BBR_{+}\ti\BBR_{+}:-\gamma x+s^{\gb}x^{p}=0\right\}=\left\{x=\left(\gamma s^{-\gb}\right)^{1/(p-1)}\right\}.
$$
For fixed $s>0$ the function $x\mapsto -2^{-1}\gamma x^2+(p+1)^{-1}s^{\gb}x^{p+1}$ vanishes at $x=0$. It has the following properties:\smallskip

\noindent (i) it is decreasing for $0<x<\left(\gamma s^{-\gb}\right)^{1/(p-1)}$,\smallskip

\noindent (ii)  it achieves a minimum at $x_{s}=\left(\gamma s^{-\gb}\right)^{1/(p-1)}$, \smallskip

\noindent (iii) and it is increasing for $x>\left(\gamma
s^{-\gb}\right)^{1/(p-1)}$ with infinite limit. Furthermore it
vanishes at $\tilde x_{s}=\left(2^{-1}(p+1)\gamma
s^{-\gb}\right)^{1/(p-1)}$.\smallskip

Let $v$ be a radially symmetric positive function. By approximation of radial elements in $H^1_{0, K}(\BBR^N)\cap L_{\abs\eta^\gb K}^{p+1}(\BBR^N)$, we can assume that $v$ is $C^2$ with nondegenerate isolated extrema. We can also assume that the graph of $v$ has at most a countable of intersections with $C_{2}$, $a_{1}<a_{2}<a_{3}...<a_{k}<...$, that the set of points $\{a_{k}\}$ is discrete, that all the intersections are transverse and that, for every $ j\geq 0$,
$$v(s)<\left(\gamma s^{-\gb}\right)^{1/(p-1)}\quad\text{on }(a_{2j},a_{2j+1}),
$$
where $a_{0}=0$, and
$$v(s)>\left(\gamma s^{-\gb}\right)^{1/(p-1)}\quad\text{on }(a_{2j+1},a_{2j+2+1}).
$$
The modifications of the function $v$ is performed by local modification on each interval $(a_{k},a_{k+1})$:\smallskip

\noindent{\it Step 1}- The construction of $\tilde v$ on $(a_{2j},a_{2j+1})$ is as follows. Let $\ga_{1}<\ga_{2}<...$ be the sequence of local extrema of $v$, with $v(\ga_{2i+1})$ local minimum and $v(\ga_{2i+2})$ local maximum. By extension, since $v'(a_{2j+1})>-\gb/(p-1)\gamma^{1/(p-1)}a_{2j+1}^{-(\gb+p-1)/(p-1)}$, $v(a_{2j+1})$ is a local maximum of $v$ on $(a_{2j},a_{2j+1})$.\smallskip

\noindent If
$\max \{(\ga_{2i+1}):i\geq 1\}\leq v(a_{2j+1})$, then $\tilde v=\max\{v,v(a_{2j+1})\}$.\smallskip

\noindent If $\max \{v(\ga_{2i+1}):i\geq 1\}> v(a_{2j+1})$, we define the increasing sequence $\{\ga_{2i_{d}+1}\}$  by
$$v(\ga_{2i_{0}+1})=\max \{v(\ga_{2i+1}):i\geq 1\},$$
$$v(\ga_{2i_{1}+1})\max \{v(\ga_{2i+1}):i> i_{0}\},$$
and by induction,
$$v(\ga_{2i_{d}+1})\max \{v(\ga_{2i+1}):i> i_{d-1}\}.$$
Thus we can assume that the local maxima of $v$ are less than $v(a_{2j+1})$ on the last interval $(\ga_{2i_{d}+1},a_{2j+1})$.
Next we define the function $\tilde v$ by   $\tilde v=\max\{v,v(\ga_{2i_{0}+1}\}$ on $(a_{2j},\ga_{2i_{0}+1})$, $\tilde v=\max\{v,v(\ga_{2i_{1}+1}\}$ on $(\ga_{2i_{0}+1},\ga_{2i_{1}+1})$. By induction, $\tilde v=\max\{v,v(\ga_{2i_{d-1}+1}\}$ on $(\ga_{2i_{d-1}+1},\ga_{2i_{d}+1})$. Finally $\tilde v=\max\{v,v(a_{2j+1})\}$ on the last interval $(\ga_{2i_{d}+1},a_{2j+1})$. The function $\tilde v$ is Lipschitz continuous, nonincreasing and, because
$v(s)\leq \tilde v(s)\leq \left(\gamma s^{-\gb}\right)^{1/(p-1)}$, there holds
\begin{equation}\label {decreas1}\BA {l}
\myint{a_{2j}\leq|\eta|\leq a_{2j+1}}{}\left(\abs{\nabla \tilde v}^2
-\gamma {\tilde v}^2+\myfrac{2}{p+1}\abs \eta^{\gb}\abs {\tilde v}^{p+1}\right)K(\eta)d\eta\\[4mm]
\qquad\quad\leq \myint{a_{2j}\leq |\eta|\leq a_{2j+1}}{}\left(\abs{\nabla v}^2
-\gamma v^2+\myfrac{2}{p+1}\abs \eta^{\gb}\abs v^{p+1}\right)K(\eta)d\eta.
\EA \end{equation}\smallskip

\noindent{\it Step 2}- The construction of $\tilde v$ on $(a_{2j+1},a_{2j+2})$ follows the same principle. Let $\gb_{1}<\gb_{2}<...<\gb_{d}$ be the sequence of local minima of $v$ on this interval. Furthermore $v(a_{2j+1})$ is the minimum of $v$ on $(a_{2j+1},a_{2j+2})$ and $v'(a_{2j+2})<-\gb/(p-1)\gamma^{1/(p-1)}a_{2j+2}^{-(\gb+p-1)/(p-1)}$.\smallskip

\noindent On $(a_{2j+1},\gb_{1})$ we set $\tilde v=\min\{v,v(a_{2j+1})\}$. On
$(\gb_{1},\gb_{2})$, $\tilde v=\min\{v,\tilde v(\gb_{1})\}$. By induction
$\tilde v=\min\{v,\tilde v(\gb_{i})\}$ on $(\gb_{i},\gb_{i+1})$. On the last interval $(\gb_{d},b_{2j+2})$, $\tilde v=\min\{v,\tilde v(\gb_{d})\}$. Because $\tilde v\leq v$ on this interval and $x\mapsto  -2^{-1}\gamma x^2+(p+1)^{-1}s^{\gb}x^{p+1}$ is increasing above the curve $C_{2}$, we obtain similarly
\begin{equation}\label {decreas2}\BA {l}
\myint{a_{2j+1}\leq |\eta|\leq a_{2j+2}}{}\left(\abs{\nabla \tilde v}^2
-\gamma {\tilde v}^2+\myfrac{2}{p+1}\abs \eta^{\gb}\abs {\tilde v}^{p+1}\right)K(\eta)d\eta\\[4mm]
\qquad\quad\leq \myint{a_{2j+1}\leq |\eta|\leq a_{2j+2}}{}\left(\abs{\nabla v}^2
-\gamma v^2+\myfrac{2}{p+1}\abs \eta^{\gb}\abs v^{p+1}\right)K(\eta)d\eta.
\EA \end{equation}\smallskip

By construction $\tilde v$ is nonincreasing. Combining (\ref{decreas1}) and (\ref{decreas2}), we obtain $J(\tilde v)\leq J(\tilde v)$.\qeda\medskip


\noindent{\it Proof of \rprop{eke}. } It follows from the previous lemmas that $J$ is bounded from below on $X$ and the function $\gf=K^{-1}$ belongs to $X$. Furthermore
$$J(t\gf)=\myfrac{(N-2\gamma)t^2}{4}\myint{}{}K^{-1}(\eta) d\eta
+\myfrac{\abs t^{p+1}}{p+1}\myint{}{}\gf^{p}(\eta) d\eta.
$$
Since $\gb>N(p-1)-2\Longleftrightarrow N-2\gamma<0$, the infimum $m$ of $J$ over radially symmetric functions is negative but finite and achieved by a decreasing function. Let
$\{v_{n}\}\subset X$ a sequence such that $J(v_{n})\downarrow m$. Then
$\{v_{n}\}$ remains bounded in $H^1_{K}(\BBR^N)\cap L_{\abs\eta^\gb K}^{p+1}(\BBR^N)$. Up to a subsequence we can assume that $v_{n}$ converges
weakly in $H^1_{K}(\BBR^N)$ and in $L_{\abs\eta^\gb K}^{p+1}(\BBR^N)$ and strongly in $L^1_{K}(\BBR^N)$ to some function $v$. Moreover this convergence holds a.e., and, since $v_{n}\in X$ the same holds with $v$. Going to the limit in the functional yields to
$$J(v)\leq\liminf_{n\to\infty}J(v_{n})=m;
$$
thus $v$ is a critical point.\qeda\medskip

The following uniqueness result holds.

\bprop {Uniq}Assume $p>1$ and $\gb> N(p-1)-2$. Then $f_{\infty}=\tilde f_{\infty}$. Furthermore $f_{\infty}$
is the unique positive solution of (\ref{self}).
\es
\Proof We first prove that $\tilde f_{\infty}$ is the unique positive radial solution of (\ref{J'})
 belonging to $H^1_{K}(\BBR^N)\cap L^{p+1}_{|\eta|^\gb K}(\BBR^N)$.
 We denote $r=|\eta|$ and $\tilde f_{\infty}(\eta)=\tilde f_{\infty}(r)$. Let $\hat f$
 be another solution in the same class. Thus there exists $\{r_{n}\}$ converging to $\infty$ such
 that $\hat f(r_{n})\to 0$. For $\ge>0$, set $\tilde f_{\ge}=\tilde f_{\infty}+\ge$. For $n\geq n_{0}$,
 large enough,
$w_{+}(r_{n})=0$, thus, as in the proof of \rprop{auxpr},
$$\BA {l}
\dint_{B_{r_{n}}}\left(\abs{\nabla \hat f-\myfrac{\hat f}{\tilde f_{\ge}}\nabla \tilde f_{\ge}}^2+\abs{\nabla \tilde f_{\ge}-\myfrac{\tilde f_{\ge}}{\hat f}\nabla {\hat f}}^2\right)Kd\eta
+\gamma\dint_{B_{r_{n}}}\myfrac{\ge}{\tilde f_{\ge}}(\hat f^{2}-\tilde f_{\ge}^{2})_{+}Kd\eta\\
\phantom{--------------------}
+\dint_{B_{r_{n}}}\abs\eta^\gb (\hat f^{p-1}-\tilde f_{\ge}^{p-1})
(\hat f^{2}-\tilde f_{\ge}^{2})_{+}Kd\eta\leq 0.
\EA$$
We let successively $r_{n}\to\infty$ with Fatou's lemma, and $\ge\to 0$ with
Lebesgue's theorem, since $\ge/\tilde f_{\ge}\leq 1$ and
$(\hat f^{2}-\tilde f_{\ge}^{2})_{+}\leq \hat f^{2}+\tilde f_{\infty}^{2}\in L^1_{K}(\BBR^N)$. We get
$$\BA {l}
\dint_{\BBR^N}\left(\abs{\nabla \hat f-\myfrac{\hat f}{\tilde f_{\infty}}\nabla \tilde f_{\infty}}^2+\abs{\nabla \tilde f_{\infty}-\myfrac{\tilde f_{\infty}}{\hat f}\nabla {\hat f}}^2+\abs\eta^\gb (\hat f^{p-1}-\tilde f_{\infty}^{p-1})
(\hat f^{2}-\tilde f_{\infty}^{2})_{+}\right)Kd\eta\leq 0,
\EA$$
which implies $\hat f\leq \tilde f_{\infty}$. In the same way
 $\tilde f_{\infty}\leq \hat f$. By \rlemma {asymp}, $f_{\infty}\in H^1_{K}(\BBR^N)\cap L^{p+1}_{|\eta|^\gb K}(\BBR^N)$. Thus $f_{\infty}=\tilde f_{\infty}$.\qeda\medskip

 We end this section with a classification result
 \bth{class} Assume $p>1$ and $\gb>N(p-1)-2$ and let $u$ be a positive solution of (\ref{E2}) which satisfies  (\ref{pc4}). Then, \smallskip

 \noindent (i) either there exists $k\geq 0$ such that $u=u_{k}$,\smallskip

 \noindent (i) or $u=u_{\infty}$.
 \es
 \Proof Because of (\ref{pc4}), the initial trace $tr(u)$ of  $u$ is is a outer regular Borel measure concentrated at $0$ (see\cite {MV2}). Then either the initial trace is a Radon measure, say $k\gd_{0}$, and we get (i), or
 \begin{equation}\label{I1}
\lim_{t\to 0}\myint{B_{\ge}}{}u(x,t)dx=\infty,
\end{equation}
for every $\ge>0$. This implies $u\geq u_{\infty}$ as in \cite {MV1}. Notice that, in this article, this estimate is performed in the case $\gb=0$, but the proof in the general case is the same. In order to prove that $u\leq u_{\infty}$, we consider, for $\ge>0$, the minimal solution $v:=v_{\ge}$ of
 \begin{equation}\label{I2}\left\{\BA {l}
\prt_{t}v-\Gd v+|x|^{\gb}|v|^{p-1}v=0\quad\text{in }\,Q_{T}\\
\phantom{------,,,--}\!\!
tr(v)=\gn_{\bar B_{\ge}},
\EA\right.\end{equation}
where $\gn_{\bar B_{\ge}}$ is the outer regular Borel measure such that
$\gn_{\bar B_{\ge}}(E)=0$ for any Borel set $E\subset\BBR^N$ such that
$E\cap \bar B_{\ge}=\emptyset$, and $\gn_{\bar B_{\ge}}(E)=\infty$ otherwhile. This solution is constructed as the limit, when $m\to\infty$ of the solution $v_{\ge,m}$ of (\ref{E2}) verifying $v_{\ge,m}(.,0)=m\chi_{_{\bar B_{\ge}}}$. Clearly $u\leq v_{\ge}$. Furthermore, for any $\ell>0$,
 \begin{equation}\label{I3}
 T_{\ell}[v_{\ge,m}]=v_{\ge/\sqrt\ell,m\ell^{(2+\gb)/2(p-1)}}
 \Longrightarrow T_{\ell}[v_{\ge}]=v_{\ge/\sqrt\ell}\Longrightarrow
 T_{\ell}[v_{0}]=v_{0},
\end{equation}
where $v_{0}=\lim_{\ge\to 0}v_{\ge}$. This, and the fact that
$\lim_{t\to 0}v_{0}(x,t)=0$ for every $x\in\BBR^N\setminus\{0\}$,
imply that $v_{0}(x,t)=t^{-(2+\gb)/2(p-1)}f_{\infty}(x/\sqrt
t)=u_{\infty}(x,t)$. At the end,  since $u\leq v_{\ge}\Longrightarrow
u\leq v_{0}$, it follows $u\leq u_{\infty}$.\qeda

\section {Existence of very singular solutions} \setcounter{equation}{0}
Our study of the singularity set of the solution $u_\infty$ in the
case of strongly degenerate potential \eqref{Flat} is based on some
variant of the local energy estimate (abr. L.E.E.) method.
First the L.E.E. method for the study of singular solutions of quasilinear
parabolic equations was used in \cite{Sh1}. Adaption of this
method to the study of conditions of removability of the point
singularities of solutions of the quasilinear parabolic equations
of diffusion-strong absorption type was given in \cite{GSh}. In
 \cite{ShV} there was elaborated a variant of the L.E.E. method, which
allowed to find sharp conditions on the time dependent absorption
potential, guaranteing existence of very singular solutions of the
Cauchy problem to diffusion-strong absorption type equation with
point singularity set. Here we provide a new application of the L.E.E. method
in describing the transformation of V.S.S solution into the R.B.
solution in terms of the flatness of the
absorption potential in the space variables.

We consider the sequence of the Cauchy problems
\begin{equation}\label{eq4.1}
u_t-\Delta u+h(|x|)|u|^{p-1}u=0 \ \ \ \text{in} \ \ \
\mathbb{R}^N\times(0,T), \ p>1,
\end{equation}
\begin{equation}\label{eq4.2}
u(x,0)=u_{0,k}(x)=M_k\,\exp (-2^{-1}\mu_0 N k) \delta_k(x),
\end{equation}
where $\delta_k$ is a regularized Dirac measure: $\delta_k\in
C(\mathbb{R}^N), \ \delta_k\rightharpoonup\delta$ weakly in the
sense of measures as $k\to\infty$,
\begin{equation}\label{eq4.3}
\text{supp}\,\delta_k\subset\{x:|x|\leq\exp (-\mu_0k)\} \ \ \ \forall
k\in\mathbb{N},
\end{equation}
where the constant $\mu_0>0$ will be defined later on, and
\begin{equation}\label{eq4.4}
M_k=\exp\,\exp k \ \ \ \forall k\in\mathbb{N}.
\end{equation}
Without loss of generality we suppose that
\begin{equation}\label{eq4.5}
\|\delta_k\|^2_{L_2(\mathbb{R}^N)}\leq \exp (\mu_0 Nk).
\end{equation}
We write the potential $h$ in the equation \eqref{eq4.1} under the
form,
\begin{equation}\label{eq4.6}
h(s)=\exp(-\omega(s)s^{-2}) \ \ \ \forall s\geq 0,
\end{equation}
where $\omega(s)\geq 0$ is arbitrary nondecreasing function on
$[0,\infty)$.

\bth{t:4.1} Let the function $\omega(s)$ defined in \eqref{eq4.6} satisfy
additionally the following Dini-like condition
\begin{equation}\label{eq4.7}
\int^{d_1}_0 \omega(s) s^{-1}ds \leq d_2<\infty, \ \ \
 d_1=\rm{const} >0,
 \end{equation}
 and the following technical condition
\begin{equation}\label{tech}
\frac{s \omega'(s)}{\omega(s)}\leq 2-\alpha_0 \ \ \ \forall
s\in(0,s_0), \ s_0>0, \ 0<\alpha_0=\rm{const}<2
 \end{equation}
 Then the following a
 priori estimate of solutions $u_k$ of the problem
 \eqref{eq4.1}, \eqref{eq4.2}, \eqref{eq4.5},  holds uniformly with respect to $k\in\mathbb{N}$,
\begin{equation}\label{eq4.8}
\int_{\mathbb{R}^N}|u_k(x,t)|^2 dx \leq C_1 t \exp
\left[C_2\left(\Phi^{-1}\left(\frac{t}{C_3}\right)\right)^{-2}\right],
\end{equation}
where the constants $C_1>0, \, C_2>0,\, C_3>0$ do not depend on $k$. Here
$\Phi^{-1}(s)$ is the inverse function to
$$
s\mapsto\Phi(s) : = \int^s_0 \frac{\omega(r)}{r} d\tau.
$$
 \es

Let us define the
 following families of domains
$$
B(s) : = \{x : |x|<s\}, \ \ \Omega(s) : = \mathbb R^N\setminus
B(s),
$$
$$
Q^{t_2}_{t_1}(s) : = \Omega(s) \times (t_1, t_2), \ \ \forall s>0,
\ \ \forall 0\leq t_1<t_2\leq T.
$$
Let $u(x,t)\equiv u_k(x,t)$ be a solution of the problem
\eqref{eq4.1}, \eqref{eq4.2} under consideration. We introduce the energy
functions
\begin{equation}\label{eq4.9*}
I(s,\tau) : = \int^\tau_0 \int_{\Omega(s)} \left(|\nabla_x u|^2 +
h(|x|) |u|^{p+1}\right) dx\,dt,
\end{equation}
\begin{equation}\label{eq4.9}
J(s,t) = \int_{\Omega(s)} |u(x,t)|^2 dx, \ \ \ E(s,t) =
\int_{B(s)} |u(x,t)|^2 dx.
\end{equation}

\blemma {l:4.1} The energy functions $J(s,t), \, I(s,t)$ defined by
\eqref{eq4.9*}, \eqref{eq4.9} corresponding to an arbitrary solution $u=u_k$ of  problem \eqref{eq4.1}, \eqref{eq4.2} satisfy the following
{\it a priori} estimate
\begin{equation}\label{eq4.10}
J(s,t)+I(s,t)\leq c t g (s) : = c t \left(\int^s_0
r^{-\frac{(N-1)(p-1)}{p+3}} h(r)^{\frac{2}{p+3}}
dr\right)^{-\frac{p+3}{p-1}}, \ \ \ \forall s \geq \exp (-\mu_0
k).
\end{equation}
uniformly with respect to $k\in\mathbb{N}$.
\es

\noindent By $c, c_i$ we denote different
positive constants, which depend on known parameters $N, p,
\alpha_0, d_2$ only, but their value may change from lines to lines.\medskip

\noindent \Proof Multiplying equation \eqref{eq4.1} by $u$ and
integrating in $Q^{t_2}_{t_1}(s)$, we obtain the following starting relation after standard computations,
\begin{multline}\label{eq4.11}
2^{-1}\int_{\Omega(s)} |u(x,t_2)|^2 dx + \iint_{Q^{t_2}_{t_1}(s)}
\left(|\nabla_xu|^2+h(|x|) |u|^{p+1}\right) dx\,dt = \\
= 2^{-1} \int_{\Omega(s)} |u(x,t_1)|^2 dx +
\int^{t_2}_{t_1}\int_{|x|=s} u \frac{\partial u}{\partial
n}\, d\sigma\,dt : = R_0+R_1.
\end{multline}
Let us estimate $R_1$ from above. Using Holder's and Young's
inequalities we have
\begin{multline*}
\left|\int_{|x|=s}u(x,t) \,\frac{\partial u}{\partial
n}\,d\sigma\right| \leq c s^{\frac{(N-1)(p-1)}{2(p+1)}}
\left(\int_{|x|=s} |\nabla_x u|^2 d\sigma\right)^{1/2}
\left(\int_{|x|=\tau} |u|^{p+1} d\sigma\right)^{\frac{1}{p+1}}
\leq \\
\leq c s^{\frac{(N-1)(p-1)}{2(p+1)}}
h(s)^{-\frac{1}{p-1}}\left(\int_{|x|=s}\left( |\nabla_x u|^2 +
h(s) |u|^{p+1}\right) d\sigma\right)^{\frac{p+3}{2(p+1)}}.
\end{multline*}
Integrating in $t$, we get
\begin{multline}\label{eq4.12}
\left|\int^\tau_0 \int_{|x|=s}u \frac{\partial u}{\partial
n}\,d\sigma\,dt\right| \leq c s^{\frac{(N-1)(p-1)}{2(p+1)}}
h(s)^{-\frac{1}{p-1}} \tau^{\frac{p-1}{2(p+1)}}  \\
 \ti\left(\int^\tau_0 \int_{|x|=s}\left( |\nabla_x u|^2 + h(s)
|u|^{p+1}\right) d\sigma\,dt\right)^{\frac{p+3}{2(p+1)}}.
\end{multline}
It is easy to see that
$$
-\frac{d}{ds} \,I(s,\tau) = \int^\tau_0 \int_{|x|=s}
\left(|\nabla_xu|^2+h(s) |u|^{p+1}\right)\,ds, \ \ \
-\frac{d}{ds}\,J(s,t)\geq 0.
$$
Therefore because of the property \eqref{eq4.3} satisfied by $u_{0,k}$, and estimate \eqref{eq4.12}, we derive the following inequality  from relation \eqref{eq4.11} with $t_2=t, \,t_1=0, \,s\geq \exp
(-\mu_0 k)$,
\begin{equation}\label{eq4.13}
J(s,t) + I(s,t) \leq c\,t^{\frac{p-1}{2(p+1)}}
h(s)^{-\frac{1}{p+1}} s^{\frac{(N-1)(p-1)}{2(p+1)}}
\left(-\frac{d}{ds} (I(s,t)+J(s,t))\right)^{\frac{p+3}{2(p+1)}}.
\end{equation}
Solving this ordinary differential inequality (abr. O.D.I.) with respect to the function $I(s,t)+J(s,t)$, we
deduce that estimate \eqref{eq4.10} holds for arbitrary $s\geq \exp
(-\mu_0k)$.  \qeda\medskip

Next, we define $s_k>0$ by the relation
\begin{equation}\label{eq4.14}
g(s_k) = M_k^{\varepsilon_0} = \exp (\varepsilon_0 \exp\,k),
\end{equation}
where $0<\varepsilon_0<1$ will be defined later on. Now we have to
guarantee that
\begin{equation}\label{Rest}
s_k\geq \exp (-\mu_0 k) : = \overline s_k \ \ \ \forall k >
k_0(\varepsilon_0, \alpha_0, \nu_0, p).
\end{equation}
Using \cite [Lemma A1]{B-Sh}, it follows from the definitions
\eqref{eq4.6}of function $h(.)$ and \eqref{eq4.10} of function
$g(.)$, that the next estimate holds,
\begin{equation}\label{eq4.15}
\left(\frac{2\alpha_0}{p+3}\right)^{\frac{p+3}{p-1}} g_1(s) \leq
g(s) \leq \left(\frac{4}{p+3}\right)^{\frac{p+3}{p-1}} g_1(s),
\end{equation}
where $g_1(s)=s^{N-1-\frac{3(p+3)}{p-1}}
\omega(s)^{\frac{p+3}{p-1}}\exp\left(\frac{2}{(p-1)}\,\frac{\omega(s)}{s^2}\right)$,
$\alpha_0$ is constant from condition \eqref{tech}. The following simpler estimate follows from \eqref{eq4.15}:
\begin{equation}\label{eq4.16}
\exp
\left(\frac{\omega(s)}{s^2}\frac{2}{(p-1)}(1-\nu_0)\right)
\leq g(s) \leq \exp\left(\frac{\omega(s)}{s^2}
\frac{2}{(p-1)}(1+\nu_0)\right),
\end{equation}
for any $s\in(0,s_0)$, where $s_0=s_0(\nu_0)\to 0$ as
$\nu_0\to 0$. As a consequence of definition \eqref{eq4.14} of
$s_k$, and using \eqref{eq4.16}, we get,
\begin{equation}\label{eq4.17}
\frac{\omega(s_k)}{s^2_k}  \frac{2(1-\nu_0)}{(p-1)} \leq
\varepsilon_0 \exp k.
\end{equation}
Integrating \eqref{tech}, we deduce that  $\omega$ satisfies
\begin{equation}\label{eq4.18}
\omega(s)\geq s^{2-\alpha_0} \ \ \ \forall s>0.
\end{equation}
Combining \eqref{eq4.18} and \eqref{eq4.17} we derive:
\begin{equation}\label{eq4.19}
s_k\geq
\left(\frac{2(1-\nu_0)}{\varepsilon_0(p-1)}\right)^{\frac{1}{\alpha_0}}
\exp \left(-\frac{k}{\alpha_0}\right).
\end{equation}
Next we define $\mu_0$ from \eqref{eq4.2} and set
$\mu_0=2\alpha_0^{-1}.$
It follows from \eqref{eq4.19} that \eqref{Rest} is
satisfied for all $k>k_0=k_0(\varepsilon_0, \alpha_0, \nu_0, p)$.
As result we derive that estimate \eqref{eq4.10} obtained in
\rlemma{l:4.1} is valid for $s=s_k$, i.e.
\begin{equation}\label{eq4.21}
J(s_k,t)+I(s_k,t)\leq ctg(s_k) \ \ \ \forall k\geq
k_0=k_0(\varepsilon_0,\alpha_0,\nu_0,p).
\end{equation}

In order to find estimates characterizing the
behaviour of the energy function $E(s_k,t)$ with respect to the variable
$t>0$, we introduce  the nonnegative cut-off function $\varphi_k\in C^1(\mathbb R)$ defined by
\begin{equation}\label{eq4.22}
\varphi_k(s)=1 \ \ \text{if} \ \ s<s_k, \ \ \ \varphi_k(s)=0 \ \
\text{if} \ \ s\geq 2s_k, \ \varphi'_k(s)\leq cs_k^{-1}.
\end{equation}
Multiplying \eqref{eq4.1} by $u_k\varphi_k^2(|x|)$
and integrating with respect to $x$, we get
\begin{multline}\label{eq4.23}
2^{-1} \frac{d}{dt} \int_{\mathbb R^N} u^2(x,t) \varphi^2_k(|x|)
dx + \int_{\mathbb R^N} |\nabla_x(u\varphi_k)|^2 dx +
\int_{\mathbb R^N} h(|x|) \varphi^2_k |u|^{p+1} dx \\
\leq \int_{\mathbb R^N} u^2(x,t)  |\nabla_x\varphi_k(|x|)|^2
dx : = \mathbb R_1.
\end{multline}
By \eqref{eq4.22} and \eqref{eq4.21}, we obtain
\begin{equation}\label{eq4.24}
\mathbb R_1 \leq c_1 s_k^{-2} \int_{s_k< |x|<2s_k} |u(x,t)|^2 dx
\leq c_1 s_k^{-2} J(s_k,t) \leq c_2 s_k^{-2} tg(s_k).
\end{equation}
Using \eqref{eq4.23}, \eqref{eq4.24} and
Poincar\'e's inequality we derive the following differential inequality,
\begin{equation}\label{eq4.25}
\frac{d}{dt} \left(\int_{\mathbb R^N} u^2(x,t) \varphi^2_k
dx\right) + d_0 s_k^{-2} \int_{B(2s_k)} u^2(x,t) \varphi^2_k dx
\leq \overline c s_k^{-2} tg(s_k), \ \ d_0>0.
\end{equation}
We set
$$\psi_k(t) : = \int_{\mathbb R^N}
|u_k(x,t)|^2 \varphi^2_k(|x|) dx,$$
and obtain the following O.D.I. from \eqref{eq4.25},
\begin{equation}\label{eq4.26}
\psi'_k(t) + d_0s_k^{-2} \psi_k(t) \leq \overline c s_k^{-2} tg
(s_k).
\end{equation}
We rewrite \eqref{eq4.26} under the form
\begin{equation}\label{eq4.27}
\psi'_k(t) + \frac{d_0}{2} s_k^{-2} \psi_k(t) + 2^{-1} \left(d_0
s_k^{-2} \psi_k(t) -2\overline c s_k^{-2} tg(s_k)\right) \leq 0.
\end{equation}
Using the relations \eqref{eq4.2}, \eqref{eq4.5} satisfied by
$u_{k,0}$, we see that $\psi_k$ verifies,
\begin{equation}\label{eq4.28}
\psi_k(0)\leq \int_{\mathbb R^N} |u_{k,0}(x)|^2 dx \leq M_k.
\end{equation}
At last, we define the $t_k$ by
\begin{equation}\label{eq4.29}
t_k=\gamma\omega(s_k)
\end{equation}
where $\omega$ is the function in \eqref{eq4.6} and $\gamma>0$ is
a parameter which will be made precise in the next lemma.

\blemma{l:4.2} There exists a constant $\gamma>0$, which does not
depend on $k$, such that any
solution $\psi_k$ of problem \eqref{eq4.27}, \eqref{eq4.28}
satisfies the following {\it a priori} estimate
\begin{equation}\label{eq4.30}
\psi_k(\overline t_k)\leq 2d_0^{-1}\overline c\, \overline t_k
g(s_k) \ \ \ \forall_k>\overline k(\varepsilon_0,\nu_0),
\end{equation}
for some $\overline t_k\leq t_k$, where $t_k$ is defined by \eqref{eq4.29}.
 \es

\noindent\Proof Let us assume that \eqref{eq4.30} is not true, and for any
$\gamma>0$ there exist $k\geq k_0$ such that
\begin{equation}\label{eq4.31}
\psi_k(t) > 2d_0^{-1}\overline c tg(s_k) \ \ \forall t :
0<t<\gamma\omega(s_k)\equiv t_k.
\end{equation}
This relation combined with \eqref{eq4.27} implies the following inequality,
$$
\psi'_k(t) +\frac{d_0}{2} s_k^{-2} \psi_k(t)\leq 0 \ \ \ \forall t
: 0<t\leq\gamma\omega(s_k).
$$
Solving this O.D.I. and using \eqref{eq4.28}, we
get
\begin{equation}\label{eq4.32}
\psi_k(t)\leq \psi_k(0) \exp \left(-\frac{d_0t}{2s_k^2}\right)
\leq M_k \exp \left(-\frac{d_0t}{2s_k^2}\right) \ \  \forall t
\leq \gamma\omega(s_k).
\end{equation}
We derive easily the next estimate from \eqref{eq4.32} and \eqref{eq4.31}
\begin{equation}\label{eq4.33}
M_k \exp \left(\frac{-d_0\gamma\omega(s_k)}{2s_k^2}\right) \geq
2d_0^{-1} \overline c g(s_k)\gamma\omega(s_k).
\end{equation}
Using \eqref{eq4.14} and \eqref{eq4.4}, we deduce from this last inequality,
\begin{equation}\label{eq4.34}
(1-\varepsilon_0) \exp k \geq \frac{d_0\gamma\omega(s_k)}{2s_k^2}
+ \ln (2d_0^{-1}\overline c\gamma) - \ln(\omega (s_k)^{-1}).
\end{equation}
Similarly to \eqref{eq4.17}, it follows, from
\eqref{eq4.16} and the definition \eqref{eq4.14} of $s_k$, that there holds
\begin{equation}\label{eq4.35}
\frac{\omega(s_k)}{s_k^2}
\frac{2(1+\nu_0)}{(p-1)}\geq\varepsilon_0 \exp k.
\end{equation}
Using this estimate and \eqref{eq4.34}, we derive
\begin{equation}\label{eq4.36} (1-\varepsilon_0) \exp k \geq
\frac{d_0\gamma (p-1)\varepsilon_0}{4(1+\nu_0)} \exp k + \ln
(d_0^{-1} 2\overline c\gamma) - \ln(\omega (s_k))^{-1}.
\end{equation}
 Noticing that \eqref{eq4.18} implies
\begin{equation}\label{eq4.37}
\ln(\omega(s_k))^{-1} \leq (2-\alpha_0) \ln(s_k^{-1}),
\end{equation}
and \eqref{eq4.19} can be writen under the form
\begin{equation}\label{eq4.38}
\ln(s_k^{-1}) \leq \frac{1}{\alpha_0}
\ln\left(\frac{\varepsilon_0(p-1)}{2(1-\nu_0)}\right) +
\frac{k}{\alpha_0},
\end{equation}
we deduce the following inequality from \eqref{eq4.37},
\eqref{eq4.38} and \eqref{eq4.36},
\begin{multline}\label{eq4.39} (1-\varepsilon_0)
\exp k \geq \frac{d_0\gamma(p-1)\varepsilon_0}{4(1+\nu_0)}
\exp k + \ln (2d_0^{-1}\overline c \gamma) - (2-\alpha_0)
\frac{k}{\alpha_0}
\\ - \frac{(2-\alpha_0)}{\alpha_0} \ln
\left(\frac{\varepsilon_0(p-1)}{2(1-\nu_0)}\right).
\end{multline}
If we define $\gamma$ by the equality
\begin{equation}\label{eq4.40}
(1-\varepsilon_0) = \frac{d_0\gamma(p-1)\varepsilon_0}{8(1+\nu_0)}
\Leftrightarrow \gamma =
\frac{(1-\varepsilon_0)(1+\nu_0)8}{d_0(p-1)\varepsilon_0} : =
\gamma_0,
\end{equation}
then inequality \eqref{eq4.39} yields to
$$
\frac{(2-\alpha_0)}{\alpha_0}k \geq (1-\varepsilon_0) \exp k + \ln
(2d_0^{-1}\overline c \gamma_0) - \frac{(2-\alpha_0)}{\alpha_0}
\ln \left(\frac{\varepsilon_0(p-1)}{2(1-\nu_0)}\right).
$$
It is clear that we can find $\overline k=\overline
k(\varepsilon_0,\nu_0)<\infty$ such that the last inequality becomes
impossible for $k\geq \overline k$, contradiction.
Consequently, \eqref{eq4.31} does not hold for
$\gamma=\gamma_0$ and estimate \eqref{eq4.30} is true with
$\gamma=\gamma_0$.  \qeda\medskip

\noindent{\it Proof of Theorem 3.1}. Comparing definition
\eqref{eq4.9} of $E(s,t)$ and definition of
$\psi_k$, we easily see that
\begin{equation}\label{eq4.41}
E(s_k,t)\leq \psi_k(t) \Rightarrow E(s_k,\overline t_k) \leq
\psi_k(\overline t_k).
\end{equation}
Therefore, using estimates \eqref{eq4.10}, \eqref{eq4.30} and
\eqref{eq4.41}, we obtain
\begin{equation}\label{eq4.42}
\int_{\mathbb R^N} |u_k(x,\overline t_k)|^2 dx = E (s_k, \overline
t_k) + J(s_k, \overline t_k) \leq
(d_0^{-1}\overline{c}+c)\overline t_k g(s_k).
\end{equation}
Next we estimate the right-hand side of \eqref{eq4.42}. Using
 \eqref{eq4.14}, \eqref{eq4.29} and inequality
\eqref{eq4.30}, we get
\begin{equation}\label{eq4.43}
\overline t_k g(s_k) \leq \gamma_0\omega(s_k) M_k^{\varepsilon_0}
\leq \gamma_0\omega(s_0)  \exp (\varepsilon_0 \exp k),
\end{equation}
where $\gamma_0$ is defined by \eqref{eq4.40}and $s_0>0$ by
\eqref{tech}. We obtain easily from \eqref{eq4.43}
\begin{equation}\label{eq4.44}
(\overline c d_0^{-1}+c)\overline t_k g(s_k) \leq \exp
\left[\left(\varepsilon_0+\frac{\ln
(\gamma_0\omega(s_0)(c+\overline{c}d_0^{-1}))}{\exp k}\right) \exp
k\right].
\end{equation}
Let  $k_1$ be the smallest integer such that
\begin{equation}\label{eq4.45}
\ln\left(\gamma_0\omega(s_0)(c+\overline{c}d_0^{-1})\right) \leq
\varepsilon_0 \exp k_1,
\end{equation}
equivalently
$$k_1 =
\left[\ln\left(\varepsilon_0^{-1}\ln\left(\gamma_0
\omega(s_0)(c+\overline{c}d_0^{-1})\right)\right)\right]+1,
$$
where $[a]$ denote integer part of $a$. Then it follows from
\eqref{eq4.44}
\begin{equation}\label{eq4.46}
(\overline c d_0^{-1}+c)t_k g(s_k) \leq \exp (2\varepsilon_0 \exp
k) \ \ \ \forall k>k_1.
\end{equation}
If we fix $\varepsilon_0$ such that
\begin{equation}\label{eq4.47}
2\varepsilon_0\leq e^{-1},
\end{equation}
then the next  estimate follows from \eqref{eq4.42} and
\eqref{eq4.43}--\eqref{eq4.47}
\begin{equation}\label{eq4.48}
\int_{\mathbb R^N} |u_k(x,\overline{t}_k)|^2 dx \leq M_{k-1},
\end{equation}
for all $k\geq \max \{k_0, \overline k, k_1\}$, where $k_0$ is
from  (\ref{Rest}), $\overline k$ -- from \eqref{eq4.30}, and $k_1$
 from \eqref{eq4.45}. Estimate \eqref{eq4.48} is the final step of
the first round of computations. For the second round, we begin by
definiting $s_{k-1}$ analogously to $s_k$:
\begin{equation}\label{eq4.49}
g(s_{k-1}) = M^{\varepsilon_0}_{k-1} = \exp (\varepsilon_0 \exp
(k-1)).
\end{equation}
From estimate \eqref{eq4.10} we obtain
\begin{equation}\label{49*}
J(s_{k-1},t)+I(s_{k-1},t)\leq ctg(s_{k-1}),
\end{equation}
since $s_{k-1}>s_k$. Analogously to
$\varphi_k$, we define the function $\varphi_{k-1}$ and set
$$
\psi_{k-1}(t) : = \int_{\mathbb R^N} |u_k(x,t)|^2
|\varphi_{k-1}(x)|^2 dx.
$$
In the same way as \eqref{eq4.26}, the following O.D.I. follows
\begin{equation}\label{eq4.50}
\psi'_{k-1}(t)+ d_0 s^{-2}_{k-1} \psi_{k-1}(t) \leq \overline c
s^{-2}_{k-1} tg(s_{k-1}) \ \ \ \forall t > \overline t_k.
\end{equation}
Using \eqref{eq4.48}, we derive
\begin{equation}\label{eq4.51}
\psi_{k-1}(\overline t_k) \leq M_{k-1}, \ \ \  \overline t_k\leq
t_k.
\end{equation}
If we analyze the Cauchy problem \eqref{eq4.50}, \eqref{eq4.51} similarly
as  problem \eqref{eq4.26}, \eqref{eq4.28}) was analyzed in
\rlemma{l:4.2},
 we obtain the following {\it a priori} estimate for $\psi_{k-1}(t)$,
\begin{equation}\label{eq4.52}
\psi_{k-1}(\overline t_k+\overline t_{k-1}) \leq
2d_0^{-1}\overline c (\overline t_k+\overline t_{k-1}) g(s_{k-1}),
\end{equation}
where $\overline t_{k-1} \leq t_{k-1} : = \gamma_0\omega(s_{k-1}),
\ \gamma_0$ is from \eqref{eq4.40}. It is clear that
$$
E(s_{k-1},t)\leq \psi_{k-1}(t) \ \ \ \forall t\geq \overline t_k,
$$
consequently
\begin{equation}\label{eq4.53}
E(s_{k-1}, \overline t_k+\overline t_{k-1}) \leq
\psi_{k-1}(\overline t_k+\overline t_{k-1}) \leq
2d_0^{-1}\overline c (\overline t_k+\overline t_{k-1}) g(s_{k-1}).
\end{equation}
From (\ref{49*}), we deduce
\begin{equation}\label{eq4.54}
J(s_{k-1}, \overline t_k+\overline t_{k-1}) + I(s_{k-1}, \overline
t_k+\overline t_{k-1}) \leq c (\overline t_k+\overline t_{k-1})
g(s_{k-1}).
\end{equation}
Summing estimates \eqref{eq4.53} and \eqref{eq4.54} we obtain
\begin{equation}\label{eq4.55}
\int_{\mathbb R^N} |u_k(x,\overline t_k+\overline t_{k-1})|^2 dx
\leq (\overline c d_0^{-1}+c) (\overline t_k+\overline t_{k-1})
g(s_{k-1}),
\end{equation}
 and we  use this last estimate for performing a similar third round of
computations. Iterating this process $j$ times, we deduce
\begin{equation}\label{eq4.56}
\int_{\mathbb R^N} \left|u_k\left(x, \sum^{k-j}_{i=k} \overline
t_i\right)\right|^2 dx \leq (\overline c d_0^{-1}+c)
\left(\sum^{k-j}_{i=k} \overline t_i\right) g(s_{k-j}).
\end{equation}
In particular, we can take  $j=k-l$, where $l\in N$ satisfies

\begin{equation}\label{eq4.57}
l\geq l_0 : = \max \{k_0, \overline k, k_1\}.
\end{equation}
Then we obtain:
\begin{equation}\label{eq4.58}
\int_{\mathbb R^N} \left|u_k\left(x, \sum^{l}_{i=k} \overline
t_i\right)\right|^2 dx \leq (\overline c d_0^{-1}+c)
\left(\sum^{l}_{i=k} \overline t_i\right) g(s_{l}).
\end{equation}
Next, we have to estimate from above the sum of the $\overline t_i$ for which there holds
\begin{equation}\label{eq4.59}
\sum^l_{i=k} \overline t_i \leq \sum^l_{i=k} \gamma_0 \omega(s_i),
\end{equation}
where $s_i$ is defined by
$g(s_i) = M_i^{\varepsilon_0}.$
By the same way as in \eqref{eq4.35}, we obtain
$$
s^2_i \leq \frac{2(1+\nu_0)\omega(s_i)}{(p-1)\varepsilon_0}
\exp(-i)\leq \frac{2(1+\nu_0)\omega(s_0)}{(p-1)\varepsilon_0}
\exp(-i) \ \ \ \forall i\geq l_0,
$$
 where $l_0$ is the integer appearing in (\ref{eq4.57}), and from this inequality follows
\begin{equation}\label{eq4.60}
s_i\leq
\left(\frac{2(1+\nu_0)\omega(s_0)}{(p-1)\varepsilon_0}\right)^{1/2}
\exp\left(-\frac{i}{2}\right) : =
C_1\exp\left(-\frac{i}{2}\right).
\end{equation}
Therefore, using the monotonicity of the function $\omega$, we derive
\begin{equation}\label{eq4.61}\BA {l}\displaystyle
\sum^l_{i=k} \omega(s_i) \leq\sum^l_{i=k} \omega
\left(C_1\exp\left(-\frac{i}{2}\right)\right)
\leq
-\int^{l-1}_{k}\omega
\left(C_1\exp\left(-\frac{s}{2}\right)\right) ds \\[2mm]
\displaystyle\phantom{\sum^l_{i=k} \omega(s_i) \leq\sum^l_{i=k} \omega
\left(C_1\exp\left(-\frac{i}{2}\right)\right)}
\leq
2 \int^{C_1\exp(-\frac{l-1}{2})}_{C_1\exp(-\frac{k}{2})}
y^{-1}\omega(y) dy\\[2mm]
 \displaystyle\phantom{\sum^l_{i=k} \omega(s_i) \leq\sum^l_{i=k} \omega
\left(C_1\exp\left(-\frac{i}{2}\right)\right)}\leq 2\int^{C_1\exp(-\frac{l-1}{2})}_0
y^{-1}\omega(y) dy \\[2mm]
\displaystyle\phantom{\sum^l_{i=k} \omega(s_i) \leq\sum^l_{i=k} \omega
\left(C_1\exp\left(-\frac{i}{2}\right)\right)}
: = 2\Phi \left(C_1\exp(-\frac{l-1}{2})\right).
\EA\end{equation}
As a consequence  of \eqref{eq4.59} and \eqref {eq4.61}, we get
\begin{equation}\label{eq4.62}
\sum^l_{i=k}\overline t_i \leq \sum^l_{i=\infty} t_i \leq
2\gamma_0 \Phi\left(C_1\exp(-\frac{l-1}{2})\right) : = T_l.
\end{equation}
The Dini condition \eqref{eq4.7} implies that $T_l\to 0 \ \
\text{as} \ \  l\to\infty$. Next, we deduce from \eqref{eq4.58} that
\begin{equation}\label{eq4.64}
\int_{\mathbb R^N} |u_k(x,T_l)|^2 dx \leq C_2 T_l g(s_l), \ \
C_2=\overline c d_0^{-1} + c \ \ \ \forall k \geq l\geq l_0.
\end{equation}
Using the fact that $s_l : g(s_l) = M_l^{\varepsilon_0}$ and \eqref{eq4.64}, we derive
\begin{equation}\label{eq4.65}
\int_{\mathbb R^N} |u_k(x,T_l)|^2 dx \leq C_2 T_l \exp
(\varepsilon_0 \exp l).
\end{equation}
Because \eqref{eq4.62} implies
\begin{equation}\label{eq4.66}
\exp l = e C^2_1
\left(\Phi^{-1}\left(\frac{T_l}{2\gamma_0}\right)\right)^{-2},
\end{equation}
we get the following inequality by plugging last relationship
into \eqref{eq4.65}:
$$
\int_{\mathbb R^N} |u_k(x,T_l)|^2 dx \leq C_2 T_l \exp
\left[e\cdot\varepsilon_0 C^2_1
\left(\Phi^{-1}\left(\frac{T_l}{2\gamma_0}\right)\right)^{-2}\right]
\ \ \forall l \geq l_0.
$$
At last, combining last estimate with (\ref{eq4.66}), we obtain
$$
\int_{\mathbb R^N} |u_k(x,t)|^2 dx \leq C_2 t \exp
\left[e^2\cdot\varepsilon_0 C^2_1
\left(\Phi^{-1}\left(\frac{t}{2\gamma_0}\right)\right)^{-2}\right]
\ \ \forall t >0,
$$
which ends the proof.\qeda

\bex{ex:4.1} Assume $\omega(s) = s^{2-\alpha_0}, \ 0<\alpha_0<2$.
Then
$$
\Phi(s) = \int^s_0 s^{1-\alpha_0} ds =
\frac{s^{2-\alpha_0}}{2-\alpha_0} \Rightarrow \Phi^{-1}(s) =
(2-\alpha)^{\frac{1}{2-\alpha_0}} s^{\frac{1}{2-\alpha_0}}.
$$
 \es
Consequently, estimate (\ref{eq4.8}) reads as follows,
$$
\int_{\mathbb R^N} |u_k(x,t)|^2 dx \leq C_1 t \exp
\left[C_2
\left(\frac{C_3}{2-\alpha_0}\right)^{\frac{2}{2-\alpha_0}}
t^{-\frac{2}{2-\alpha_0}}\right] \ \ \forall t > 0.
$$

\section { Razor blades}
\setcounter{equation}{0}
In this section we consider equation (\ref{E1}) with potential $h(|x|)$ of the form (\ref{eq4.6})
with the limiting function $\omega(|x|):=|x|^2\ell(|x|),$ namely, we study the equation
\begin{equation}\label{F1}
\prt_tu-\Gd u+e^{-\ell(|x|)}\abs u^{p-1}u=0,\quad\text { in }\BBR^N\ti (0,\infty)
\end{equation}
where $\ell\in C(\BBR^N)$ is positive, nonincreasing and $\lim_{r\to 0}\ell(r)=\infty$. Our main result is the following
\bth{RBth} Assume $p>1$ and  $\ell$ satisfies
\begin{equation}\label{RB1}
\lim\inf_{x\to 0}|x|^2\ell(x)>0.
\end{equation}
Then the solution $u_{k}$ of the problem  (\ref{E1}), (\ref{init}),
exists for any $k>0$ and $u_{\infty}:=\lim_{k\to\infty}$ is a solution of (\ref{F1}) in $Q_{\infty}\setminus\{0\}\ti\BBR^+$ with the following properties,
\begin{equation}\label{RB2}
\lim_{t\to 0}u_{\infty}(x,t)=0\forevery x\neq 0\quad\text{and }\;\;\lim_{x\to 0}u_{\infty}u(x,t)=\infty\forevery t> 0.
\end{equation}
Furthermore $t\mapsto u_{\infty}(x,t)$ is increasing and $\lim_{t\to\infty}u_{\infty}(x,t)=U(x)$ for every $x\neq 0$ where $U=\lim_{k\to\infty}U_{k}$ and $U_{k}$ solves
\begin{equation}\label{RB3}
-\Gd U_{k}+e^{-\ell (x)}U_{k}^p=k\gd_{0}\quad\text {in }\CD'(\BBR^N).
\end{equation}
\es
\Proof By assumption (\ref{RB1}), property (\ref{int1}) is fulfilled. Thus for $k>0$ there exists
$u:=u_{k}$ solution of (\ref{F1}), (\ref{init}). Moreover, for any $k>0$ there exists
a solution $U_{k}$ of (\ref{RB3}) (see \cite{Ve}); the mapping $k\mapsto U_{k}$ is increasing and $U=\lim_{k\to\infty}U_{k}$ exists, because of Keller-Osserman estimate. $U$ is the
minimal solution of
\begin{equation}\label{RB4}
-\Gd V+e^{-\ell (x)}V^p=0\quad\text {in }\BBR^N\setminus\{0\},
\end{equation}
verifying
\begin{equation}\label{RB4'}
\myint{B_{\ge}}{}V(x)dx=\infty\forevery\ge>0.
\end{equation}
If we denote by $\bar U$ the maximal solution of (\ref{RB4}), it is classical that $\bar U=\lim_{\ge\to 0}\bar U_{\ge}$ where
\begin{equation}\label{RB4''}\left\{\BA {l}
-\Gd \bar U_{\ge}+e^{-\ell(x)}\bar U^p_{\ge}=0
\quad\text{in }\BBR^N\setminus \bar B_{\ge}\\
\phantom{--}
\lim_{|x|\to\ge}\bar U_{\ge}(x)=\infty.
\EA\right.\end{equation}
Since any $u_{k}$ is bounded from above by $\bar U$, the local equicontinuity of the $u_{k}$ in $\bar Q_{T}\setminus\{(0,0)\}$ implies that $u_{\infty}$ satisfies $\lim_{t\to 0}u_{\infty}(x,t)=0$ for all $x\neq 0$.
\smallskip

\noindent {\it Step 1: Formation of the razor blade}. \underline{The Case 1: $1<p<1+2/N$}.
For $\ge>0$, $e^{-\ell(|x|)}\leq e^{-\ell(\ge)}$ for $\abs x\leq \ge$. Therefore
\begin{equation}\label{F2}
\prt_tu-\Gd u+e^{-\ell(\ge)}\abs u^{p-1}u\geq 0,\quad\text { in }B_{\ge}\ti (0,\infty).
\end{equation}
and $u\geq v_{\ge}$ in $B_{\ge}\ti (0,T)$ where $v_{\ge}$ solves
\begin{equation}\label{F3}\left\{\BA {l}
\prt_tv_{\ge}-\Gd v_{\ge}+e^{-\ell(\ge)}\abs {v_{\ge}}^{p-1}v_{\ge}= 0\quad\text { in }B_{\ge}\ti (0,\infty)\\
\phantom{\prt_tv_{\ge}-\Gd v_{\ge}+e^{-\ell(\ge)}\abs v_{\ge}^{p-1}}
v_{\ge}=0\quad \text { in }\prt B_{\ge}\ti (0,\infty)
\\
\phantom{-;;\Gd +e^{-\ell(\ge)}\abs v_{\ge}^{p-1}}
v_{\ge}(x,0)=\infty\gd_{0}\quad \text { in }B_{\ge},
\EA\right.\end{equation}
where the initial condition is to be understood in the sense $\lim_{k\to\infty}k\gd_{0}$.
We put
$$w_{\ge}(x,t)=\ge^{2/(p-1)}e^{-\ell(\ge)/(p-1)}v_{\ge}(\ge x,\ge^2t).
$$
Then $w_{\ge}=w$ is independent of $\ge$ and solves
\begin{equation}\label{F3'}\left\{\BA {l}
\prt_tw-\Gd w+\abs w^{p-1}w= 0\quad\text { in }B_{1}\ti (0,\infty)\\
\phantom{\prt_tw-\Gd w+\abs w^{p-1}}
w=0\quad \text { in }\prt B_{1}\ti (0,\infty)
\\
\phantom{.\prt_tw+\abs w^{p-1}}
w(x,0)=\infty\gd_{0}\quad \text { in }B_{1}.
\EA\right.\end{equation}
Therefore
\begin{equation}\label{F4}u(0,1)\geq v_{\ge}(0,1)=\ge^{-2/(p-1)}e^{\ell(\ge)/(p-1)}w(0,\ge^{-2}).
\end{equation}
The longtime behaviour is given in \cite{GV} where it is proved
$$\lim_{\gt\to\infty}e^{\gl_{1}\gt}w(0,\gt)=\gk\phi_{1}(0).
$$
In this formula $\phi_{1}$ is the first eigenfunction of $-\Gd$ in $W^{1,2}_{0}(B_{1})$,  $\gl_{1}$ the corresponding eigenvalue and $\gk>0$. Thus
\begin{equation}\label{F5}
u(0,1)\geq \gd\ge^{-2/(p-1)}e^{\ell(\ge)/(p-1)}e^{\gl_{1}\ge^{-2}}\phi_{1}(0),
\end{equation}
for some $\gd>0$, if $\ge$ is small enough. If we assume
\begin{equation}\label{F6}
\lim_{\ge\to 0}\left(\myfrac{2}{p-1}\ln\ge^{-1}+\myfrac{\ell(\ge)}{p-1}-\gl_{1}\ge^{-2}\right)=\infty,
\end{equation}
it implies
\begin{equation}\label{F7}
u(0,1)=\infty\Longrightarrow u(0,t)=\infty\forevery t> 0.
\end{equation}
Moreover, the unit ball $B_{1}$ can be replaced by any ball $B_{R}$ and
$\gl_{1}$ by $\gl_{R}=R^{-2}\gl_{1}$. Therefore the sufficient condition for a Razor blade is that
it exists some $c>0$ such that
\begin{equation}\label{F8}
\lim_{\ge\to 0}\left(\ell(\ge)-c\ge^{-2}\right)=\infty.
\end{equation}
An equivalent condition is
\begin{equation}\label{F9}
\liminf_{\ge\to 0}\ge^2\ell(\ge)>0.
\end{equation}
\underline{The general case.} If $p>1$ is arbitrary, we consider $\gb>0$ such that $\gb>N(p-1)-2$, and we write
$$e^{-\ell(x)}=|x|^{\gb}e^{-\ell (x)-\gb\ln|x|}.
$$
For $R>0$ small enough $x\mapsto\tilde\ell (x):= \ell (x)+\gb\ln|x|$ is positive, increasing and satisfies the same blow-up condition (\ref{RB1}) as $\ell$. Clearly
$u_{k}$ is bounded from below on $B_{R}\ti (0,\infty)$ by the solution $\tilde u:=\tilde u_{k}$ of
\begin{equation}\label{F10}\left\{\BA {l}
\prt_t\tilde u-\Gd \tilde u+|x|^{\gb}e^{-\tilde\ell(x)}\abs {\tilde u}^{p-1}\tilde u= 0\quad\text { in }B_{R}\ti (0,\infty)\\
\phantom{|x|^{\gb}\prt_t\tilde u-\Gd \tilde u+e^{-\tilde\ell(x)}\abs {\tilde u}^{p-1}}
\tilde u=0\quad \text { in }\prt B_{R}\ti (0,\infty)
\\
\phantom{|x|^{\gb},\Gd \tilde u+e^{-\tilde\ell(x)}\abs {\tilde u}^{p-1}}
\tilde u(x,0)=k\gd_{0}\quad \text { in }B_{R}.
\EA\right.\end{equation}
Therefore, for $0<\ge<R$, $\tilde u_{\infty}$ is bounded from below on
$B_{\ge}\ti (0,\infty)$ by the solution $v_{\ge}$ of
\begin{equation}\label{11}\left\{\BA {l}
\prt_tv_{\ge}-\Gd v_{\ge}+|x|^{\gb}e^{-\tilde\ell(\ge)}\abs {v_{\ge}}^{p-1}v_{\ge}= 0\quad\text { in }B_{\ge}\ti (0,\infty)\\
\phantom{\prt_tv_{\ge}-\Gd v_{\ge}+|x|^{\gb}e^{-\tilde\ell(\ge)}\abs v_{\ge}^{p-1}}
v_{\ge}=0\quad \text { in }\prt B_{\ge}\ti (0,\infty)
\\
\phantom{-;;\Gd +|x|^{\gb}e^{-\tilde\ell(\ge)}\abs v_{\ge}^{p-1}}
v_{\ge}(x,0)=\infty\gd_{0}\quad \text { in }B_{\ge}.
\EA\right.\end{equation}
If we set
$$w_{\ge}(x,t)=\ge^{(2+\gb)/(p-1)}e^{-\ell(\ge)/(p-1)}v_{\ge}(\ge x,\ge^2t),
$$
then $w_{\ge}=w$ is independent of $\ge$ and
\begin{equation}\label{F12}\left\{\BA {l}
\prt_tw-\Gd w+|x|^{\gb}\abs w^{p-1}w= 0\quad\text { in }B_{1}\ti (0,\infty)\\
\phantom{\prt_tw-\Gd w+\abs w^{p-1}}
w=0\quad \text { in }\prt B_{1}\ti (0,\infty)
\\
\phantom{.\prt_tw+\abs w^{p-1}}
w(x,0)=\infty\gd_{0}\quad \text { in }B_{1}.
\EA\right.\end{equation}
By a straightforward adaptation of the result of \cite {GV}, there still holds
$$\lim_{\gt\to\infty}e^{\gl_{1}\gt}w(0,\gt)=\gk\phi_{1}(0)
$$
for some $\gk>0$. The remaining of the proof is the same as in case $1<p<1+2/N$.\medskip

\noindent{\it Step 2: Asymptotic behaviour.} A key observation is that, for any
$\gt>0$ and any $\ge_{0}>0$
\begin{equation}\label{F12'}\myint{\ge_{0}}{}u_{\infty}(x,\gt)dx=\infty.
\end{equation}
We give the proof in the case $1<p<1+2/N$, the general case being similar.  By step 1
\begin{equation}\label{F13}\myint{B_{\ge}}{}u(x,\gt)dx\geq
\myint{B_{\ge}}{}v_{\ge}(x,\gt)dx=\ge^{-2/(p-1)+N}e^{\ell(\ge)/(p-1)}
\myint{B_{1}}{}w(y,\ge^{-2}\gt)dy.
\end{equation}
If we fix $\gt$ and use \cite {GV}, there exists $\ge_{0}$ such that
$w(y,\ge^{-2}\gt)\geq 2^{-1}\gk e^{-\gl_{1}\ge^{-2}\gt}\gf_{1}(y)$ for $\ge\leq \ge_{0}$ and $y_{1}\in B_{1}$. Therefore
\begin{equation}\label{F14}
\myint{B_{\ge}}{}u(x,\gt)dx\geq c\ge^{-2/(p-1)+N}e^{\ell(\ge)/(p-1)-\gl_{1}\ge^{-2}\gt},
\end{equation}
for some constant $c>0$. If $\gt$ is small enough, the right-hand side of (\ref{F14}) tends to infinity as $\ge\to 0$, so does the left-hand side. This implies
(\ref{F12'}). For any $k>0$ and any $\ge>0$, there exists
$m=m(\ge)>0$ such that
$$\myint{B_{\ge}}{}\min\{u(x,\gt),m\}dx=k,
$$
thus, if we set $\gf_{m}=\min\{u(x,\gt),m\}\chi_{_{B_{\ge}}}$, then $u$ is bounded from below on $\BBR^N\ti (\gt,\infty)$ by the solution $v=v_{\ge,k}$ of
\begin{equation}\label{F15}\left\{\BA {l}
\prt_{t}v-\Gd v+e^{-\ell(x)}|v|^{p-1}v=0\quad\text {in } \BBR^N\ti (\gt,\infty)\\
\phantom{---)-----}
v(x,\gt)=\gf_{m}(x)\quad\text {in }\BBR^N.
\EA\right.\end{equation}
When $\ge\to 0$, $\gf_{m}(.)\to k\gd_{0}$ weakly in $\mathfrak M(\BBR^N)$. By standard approximation property, $v(\ge,k)\to v_{0,k}$ which is a solution of
\begin{equation}\label{F16}\left\{\BA {l}
\prt_{t}v-\Gd v+e^{-\ell(x)}|v|^{p-1}v=0\quad\text {in } \BBR^N\ti (\gt,\infty)\\
\phantom{---)-----}
v(.,\gt)=k\gd_{0}\quad\text {in }\BBR^N.
\EA\right.\end{equation}
By uniqueness, $v_{0,k}(x,t)=u_{k}(x,t-\gt)$. Letting $k\to\infty$ yields to
\begin{equation}\label{F17}
u_{\infty}(x,t+\gt)\geq u_{\infty}(x,t)\forevery (x,t)\in Q_{T}.
\end{equation}
This implies that $t\mapsto u_{\infty}(x,t)$ is increasing for every $x\in \BBR^N$. Because $u(x,t)\leq U(x)$, it is straightforward that $\lim_{x\to\infty}u(x,t)=\tilde U(x)$ exists in $\BBR^N\setminus\{0\}$.
\medskip

\noindent{\it Step 3: Identification of the limit.}
If $\gz\in C^\infty_{0}(\BBR^N\setminus\{0\}) $, there holds
$$\myint{T}{T+1}\myint{\BBR^N}{}\left(-u(x,t)\Gd\gz(x)+e^{-\ell(x)}u^p(x,t)\gz(x)\right)dx\,dt=\myint{\BBR^N}{}\left(u(x,T)-u(x,T+1)\right)\gz(x)dx.
$$
By Lebesgue's theorem
\begin{equation}\label{F18}\myint{\BBR^N}{}\left(-\tilde U(x)\Gd\gz(x)+e^{-\ell(x)}\tilde U^p(x)\gz(x)\right)dx=0,
\end{equation}
and, from (\ref{F12'}),
\begin{equation}\label{F19}\myint{\ge_{0}}{}\tilde U(x)dx=\infty,
\end{equation}
for any $\ge_{0}>0$. Therefore $\tilde U$ is a solution of the stationary equation (\ref{RB3}) in $\BBR^N\setminus\{0\}$ with a strong singularity at $0$. For $k>0$ and $\ge>0$ there exists $k(\ge)>0$ such that
$$\myint{B_{\ge}}{}U_{k(\ge)}dx=k.
$$
Let $v:=v_{k,\ge}$ be the solution of
\begin{equation}\label{F20}\left\{\BA {l}
\prt_{t}v-\Gd v+e^{-\ell(x)}|v|^{p-1}v=0\quad\text{in }\,Q_{T},\\\phantom{-----,,---}
v(.,0)=U_{k(\ge)}\chi_{_{B_{\ge}}}\quad\text{in }\,\BBR^N.
\EA\right.\end{equation}
Since $v_{k,\ge}(.,0)\leq U_{k(\ge)}(.)$,  the maximum principle implies $v_{k,\ge}\leq U_{k(\ge)}$. If we let $\ge\to 0$, $v_{k,\ge}$ converges to the solution  $u_{k}$ with initial data $k\gd_{0}$. Furthermore $k(\ge)\to\infty$ as $\ge\to 0$. Therefore
\begin{equation}\label{F21}
u_{k}(x,t)\leq  U(x)\quad\forall (x,t)\in Q_{T}.
\end{equation}
Letting successively $k\to\infty$ and $t\to\infty$ implies
\begin{equation}\label{F22}
\tilde U(x)\leq  U(x)\quad\forall x\in \BBR^N.
\end{equation}
Since $U$ is the minimal solution of (\ref{RB4}) verifying (\ref{RB4'}), it follows that $U=\tilde U$.\qeda

\end{document}